\documentclass{statsoc}

\usepackage{graphicx}
\usepackage{amssymb}
\usepackage{multirow}
\usepackage{anysize}

\newtheorem{theorem}{Theorem}
\newtheorem{prop}{Proposition}
\newtheorem{lemma}{Lemma}
\newtheorem{remark}{Remark}
\newcommand{\etal}{{\sl et al. }}

\newcommand{\var}{\mbox{var}}
\newcommand{\diag}{\mbox{diag}}
\renewcommand{\hat}{\widehat}

\newcommand{\SNR}{\mbox{SNR}}

\newcommand{\bX}{\mbox{\bf X}}
\newcommand{\bW}{\mbox{\bf W}}

\newcommand{\bQ}{\mbox{\bf Q}}

\newcommand{\bY}{\mbox{\bf Y}}

\newcommand{\bU}{\mbox{\bf U}}

\newcommand{\bA}{\mbox{\bf A}}
\newcommand{\ba}{\mbox{\bf a}}
\newcommand{\bB}{\mbox{\bf B}}
\newcommand{\bb}{\mbox{\bf b}}

\newcommand{\bD}{\mbox{\bf D}}

\newcommand{\bbeta}{\mbox{\boldmath${\beta}$}}

\newcommand{\boldeta}{\mbox{\boldmath${\eta}$}}

\newcommand{\bgamma}{\mbox{\boldmath${\gamma}$}}
\newcommand{\bSigma}{\mbox{\boldmath${\Sigma}$}}

\newcommand{\btheta}{\mbox{\boldmath${\theta}$}}

\newcommand{\sbtheta}{\mbox{\scriptsize \boldmath $\theta$}}
\newcommand{\sbeta}{\mbox{\scriptsize \boldmath $\eta$}}
\newcommand{\sbgamma}{\mbox{\scriptsize \boldmath $\gamma$}}

\newcommand{\balpha}{\mbox{\boldmath${\alpha}$}}

\newcommand{\bpi}{\mbox{\boldmath${\pi}$}}

\newcommand{\E}{\mathrm{E}}
\newcommand{\Var}{\mathrm{Var}}
\newcommand{\Cov}{\mathrm{Cov}}
\def\ANNALS{{\em Ann. Statist.}}

\def\JASA{{\em J. Am. Statist. Assoc.}}

\def\JRSSB{{\em J.Roy. Statisti. Soc. B.}}

\def\JCGS{{\em J. Comput. Graph. Stat.}}

\newcommand{\bg}{\begin{eqnarray}}
\newcommand{\ed}{\end{eqnarray}}
\newcommand{\bgn}{\begin{eqnarray*}}
\newcommand{\edn}{\end{eqnarray*}}

\title[NIS for Ultra-High Dimensional Varying Coefficient Models]{Nonparametric Independence Screening in Sparse Ultra-High Dimensional Varying Coefficient Models}
\author{Jianqing Fan}
\address{Princeton University, USA}
\author{Yunbei Ma}
\address{Southwestern University of Finance and Economics, China}
\author[J. Fan, Y. Ma and W. Dai]{and Wei Dai}
\address{Princeton University, USA}

\begin{document}

\maketitle

\begin{abstract}
The varying-coefficient model is an important nonparametric statistical model that allows us to examine how the effects of covariates vary with exposure variables.  When the number of covariates is big, the issue of variable selection arrives.   In this paper, we propose and investigate marginal nonparametric screening methods to screen variables in ultra-high dimensional sparse varying-coefficient models. The proposed nonparametric independence screening (NIS) selects variables by ranking a measure of the nonparametric marginal contributions of each covariate given the exposure variable. The sure independent screening property is established under some mild technical conditions when the dimensionality is of nonpolynomial order, and the dimensionality reduction of NIS is quantified.  To enhance practical utility and the finite sample performance,  two data-driven iterative NIS methods are proposed for selecting thresholding parameters and variables: conditional permutation and greedy methods, resulting in Conditional-INIS and Greedy-INIS. The effectiveness and flexibility of the proposed methods are further illustrated by simulation studies and real data applications.

\keywords{Sure independence screening; Variable selection; Sparsity; Conditional permutation; False positive rates}
\end{abstract}

\footnotetext
{\emph{Address for correspondence}: Yunbei Ma, School of Statistics, Southwest University of Finance and Economics, Chengdu, China.\\
E-mail:myb@swufe.edu.cn}

\section{Introduction}
The development of information and technology drives big data
collections in many areas of advanced scientific research ranging
from genomic and health science to machine learning and economics.
The collected data frequently has an ultra-high dimensionality $p$
that is allowed to diverge at nonpolynomial (NP) rate with the
sample size $n$, namely $\log(p)=O(n^\rho)$ for some $\rho>0$. For
example, in biomedical research such as genomewide association
studies for some mental diseases, millions of SNPs are potential
covariates. Traditional statistical methods face significant
challenges in dealing with such a high-dimensional problem with
large sample sizes.

With the sparsity assumption, variable selection helps improve the
accuracy of estimation and gain scientific insights. Many
significant variable selection techniques have been developed, such
as Bridge regression in \cite{Frank1993}, Lasso in
\cite{Tibshirani1996},  SCAD and folded concave penalty in
\cite{Fan2001}, the Elastic net in \cite{Zou2005}, Adaptive Lasso
\citep{Zou2006}, and the Dantzig selector in \cite{Candes2007}.
Methods on the implementation of folded concave penalized
least-squares include the local linear approximation algorithm in
\cite{Zou2008} and the plus algorithm in \cite{Zhang2010}. However,
due to the simultaneous challenges of computational expediency,
statistical accuracy and algorithmic stability, these methods do not
perform well in ultra-high dimensional problems.

To tackle these problems, \cite{Fan2008} introduced a sure
independence screening (SIS) method to select important variables in
ultra-high dimensional linear regression models via marginal
correlation learning. \cite{Hall2009} extended the method to the
generalized correlation ranking, which was further extended by
\cite{Fan2011} for ultra-high dimensional nonparametric additive
models, resulting in nonparametric independence screening (NIS).  On
a different front, \cite{Fan2010} extended the SIS idea to
ultra-high dimensional generalized linear models and devised a
useful technical tool for establishing the sure screening results
and bounding false selection rates. Other related methods include
data-tilling method \citep{Hall22009}, marginal partial likelihood
method MPLE \citep{Zhao2010}, and robust screening methods by rank
correlation \citep{Li2012} and distance correlation \citep{Li22012}.
Inspired by these previous work, our study will focus on variable
screening in nonparametric varying-coefficient models with NP
dimensionality.

It is well known that nonparametric models are flexible enough to
reduce modeing biases. However, they suffer from the so-called
``curse of dimensionality". A remarkable simple and powerful
nonparametric model for dimensionality reductions is the
varying-coefficient model,
 \bg
 Y =\bbeta^T(W)\bX + \epsilon,
 \label{eq1}
 \ed
where $\bX=(X_1,\cdots, X_p)^T$ is the vector of covariates, $W$ is
some observable exposure variables, $Y$ is the response, and
$\varepsilon$ is the random noise with conditional mean 0 and finite
conditional variance. An intercept term (i.e., $X_0 \equiv 1$) can
be introduced if necessary. This model assumes that the variables in
the covariate vector $\bX$ enter the model linearly, meanwhile it
allows regression coefficient functions to very smoothly with the
exposure variable.  The model retains general nonparametric
characteristics and allows the nonlinear interactions between the
exposure variable $W$ and the covariates.  It arises frequently from
economics, finance, politics, epidemiology, medical science,
ecology, among others.  For an overview, see \cite{Fan22008}.

When the dimensionality $p$ is finite, \cite{Fan22001} proposed the
generalized likelihood ratio (GLR) test to select variables in the
varying-coefficient model (\ref{eq1}). For the time-varying
coefficient model, a special case of (\ref{eq1}) with the exposure
variable $W$ being the time $t$, \cite{Wang2008}  applied the basis
function approximations and the SCAD penalty to address the problem
of variable selection. In the NP dimensional setting,
\cite{Lian2011} utilized the adaptive group Lasso penalty in
time-varying coefficient models. These methods still face the
aforementioned three challenges.

In this paper, we consider a nonparametric screening by ranking a measure of the marginal nonparametric contributions of each covariate given the exposure variable. For each given covariate,
we fit marginal regressions of the response $Y$ against the covariate $X_j$ $(j=1,\cdots,p)$ conditioning on $W$:
\bg
\min_{a_j, b_j} \E[(Y - a_j - b_j X_j)^2 | W]
 \label{eq2}
\ed Let $a_j(W)$ and $b_j(W)$ be the solution to (\ref{eq2}) and
$\hat{a}_{nj}(W)$ and $\hat{b}_{nj}(W)$ be their nonparametric
estimates. Then, we rank the importance of each covariate in the
joint model according to a measure of marginal utility (which is
equivalent to the goodness of fit) in its marginal model. Under some
reasonable conditions, the magnitude of these marginal contributions
provides useful probes of the importance of variables in the joint
varying-coefficient model. This is an important extension of SIS
\citep{Fan2008} to a more flexible class of varying coefficient
models.

The sure screening property of NIS can be established under certain
technical conditions. In some very specific cases, NIS can even be
model selection consistent. In establishing this kind of results,
three factors are related to the minimum distinguishable marginal
signals: the stochastic error in estimating the nonparametric
components, the approximation error in modeling nonparametric
components, and the tail distributions of the covariates. Following
\cite{Fan2008} and \cite{Fan2011}, we propose two nonparametric
independence screening approaches in an iterative framework. One is
called Greedy-INIS, in which we adopt a greedy method in the
variable screening step. The other is called Conditional-INIS which
is built on conditional random permutation to determine a data
driven screening threshold. They both serve to effectively control
the false positive rate and false negative rate with enhanced
performance.

This article is organized as follows. In Section 2, we fit each marginal nonparametric regression model via B-spline basis approximation and screen variables by ranking a measure of these estimators. In Section 3, we establish the sure screening property and model selection consistency under certain technical conditions. Iterative NIS procedures (namely Greedy-INIS and Conditional-INIS) are developed in Section 4. In Section 5, a set of numerical studies are conducted to evaluate the performance of our proposed methods.

\section{Models and Nonparametric Marginal Screening Method}
In this section we study the varying-coefficient model with the conditional linear structure as in (\ref{eq1}).
Assume that the functional coefficient vector $\bbeta(\cdot)=(\beta_1(\cdot),\cdots,\beta_p(\cdot))^T$ is sparse. Let $\mathcal{M}_*=\{j: \E[\beta^2_j(W)]>0\}$ be the true sparse model with nonsparsity size $s_n=|\mathcal{M}_*|$. We allow $p$ to grow with $n$ and denote it by $p_n$ whenever necessary.
\subsection{Marginal Regression}
For $j=1,\cdots,p$, let $a_j(W)$ and $b_j(W)$ be the minimizer of the following marginal regression problem:
 \bg
  \min\limits_{a_j(W), b_j(W) \in L_2(P)}\E[(Y- a_j(W) - b_j(W)X_j)^2|W],
  \label{eq3}
 \ed
where $P$ denotes the joint distribution of $(Y,W,\bX)$ and $L_2(P)$ is the class of square integrable functions under the measure $P$. By some algebra, we have that the minimizer of (\ref{eq3}) is
\bg
b_j(W)=\frac{\Cov[X_j,Y |W]}{\Var[X_j|W]}, ~a_j (W)= \E[Y|W] - b_j(W)\E[X_j|W].
\label{eq4}
\ed
Let $a_0(W)=\E[Y|W]$, we rank the marginal utility of covariates by
\bg
u_j = \| a_j(W) + b_j(W)X_j \|^2 - \|a_0(W)\|^2,
\label{eq5}
\ed
 where $\|f\|^2 = \E f^2.$ It can be seen that
\bg
u_j &=& \E[b_j^2(W) (X_j - \E[X_j|W])^2]
 = \E\left[ \frac{(\Cov[X_j,Y|W])^2}{\Var[X_j|W]}\right].
\label{eq6}
\ed

For each  $j=1,\cdots,p$, if $\Var[X_j |W] = 1$, then $u_j$ has the same quantity as the measure of marginal functional coefficient $\|b_j(W)\|^2$. On the other hand, this marginal utility is closely related to the conditional correlation between $X_j's$ and $Y$, as $u_j = 0 $ if and only if $\Cov[X_j,Y|W] = 0$ almost surely.

\subsection{Marginal Regression Estimation with B-spline}
\label{2.2} To obtain an estimate of the marginal utility $u_j$, $j=1, \cdots, p$, we approximate $a_j(W)$ and $b_j(W)$ by functions in
$\mathcal{S}_n$, the space of polynomial splines of degree $l \geq 1$ on $\mathcal{W}$, a compact set. Let $\{B_{k},k=1,\cdots,L_n\}$ denote its normalized B-spline basis with $\|B_{k} \|_\infty\leq 1,$ where $\|\cdot\|_\infty$ is the sup norm.  Then
\bg
a_j(W) &\approx& \sum_{k=1}^{L_n}\eta_{jk}B_{k}(W), ~\quad j=0, \cdots, p, \nonumber\\
b_j(W) &\approx& \sum_{k=1}^{L_n}\theta_{jk}B_{k}(W), ~\quad j=1, \cdots, p.\nonumber
\ed
where $\{\theta_{jk}\}_{k=1}^{L_n}$ and $\{\eta_{jk}\}_{k=1}^{L_n}$ are scalar coefficients.

We now consider the following sample version of the marginal regression problem:
 \bg
   \min\limits_{\sbeta_j,\sbtheta_j\in \mathbb{R}^{L_n}}\frac1n\sum_{i=1}^n(Y_i-\bB(W_i)\boldeta_j-\bB(W_i)\btheta_{j} X_{ji})^2,
\label{eq7}
 \ed
where $\boldeta_j=(\eta_{j1},\cdots,\eta_{jL_n})^T$,  $\btheta_j=(\theta_{j1},\cdots,\theta_{jL_n})^T$ and $\bB(\cdot)=(B_{1}(\cdot),\cdots,B_{L_n}(\cdot))$.

It is easy to show that the minimizers of (\ref{eq7}) is given by
\begin{equation}
 (\hat{\boldeta}^T_{j},~ \hat{\btheta}^T_j)^T=(\bQ^T_{nj}\bQ_{nj})^{-1}\bQ^T_{nj}\bY,
\label{eq8}
\end{equation}
where
\bg
\bQ_{nj}=\left( \bB_n,~ \bold{\Phi}_{nj}\right) = \left(\begin{array}{cc}\bB(W_1),&X_{j1}\bB(W_1)\\\vdots&\vdots\\\bB(W_n),&X_{jn}\bB(W_n)\end{array}\right) \nonumber
\ed
is an $n \times 2L_n$  matrix. As a result, the  estimates of $a_j$ and $b_j$, $j=1,\cdots,p$ are given by
\bg
 &\hat{a}_{nj}(W) = \bB(W) \hat{\boldeta}_j = (\bB(W),{\bf 0}^T_{L_n})(\bQ^T_{nj}\bQ_{nj})^{-1}\bQ^T_{nj}\bY, & \nonumber\\
 &\hat{b}_{nj}(W) = \bB(W) \hat{\btheta}_j = ({\bf 0}^T_{L_n},\bB(W))(\bQ^T_{nj}\bQ_{nj})^{-1}\bQ^T_{nj}\bY,&
\label{eq9}
\ed
where ${\bf 0}_{L_n}$ is an $L_n$-dimension vector with all entries 0. Similarly, we have the estimate of the intercept function $a_0$ by
\bg
\hat{a}_{n0}(W) = \bB(W) \hat{\boldeta}_0 = \bB(W)(\bB^T_n\bB_n)^{-1}\bB^T_n\bY,
\label{eq10}
\ed
where
\bg
\hat{\boldeta}_0 =  \arg \min\limits_{\sbeta_0 \in \mathbb{R}^{L_n}}\frac1n\sum_{i=1}^n(Y_i-\bB(W_i)\boldeta_0)^2.
\label{eq11}
\ed

We now define an estimate of the marginal utility $u_j$ as
\bg
\hat{u}_{nj} &=& \|\hat{a}_{nj}(\bW) + \hat{b}_{nj} (\bW)\bX_j \|^2_n - \| \hat{a}_{n0}(\bW) \|^2_n \nonumber \\
&=&\frac1n\sum_{i=1}^n(\hat{a}_{nj}(W_i) + \hat{b}_{nj}(W_i)X_{ji})^2 - \frac1n\sum_{i=1}^n(\hat{a}_{n0}(W_i))^2,
\label{eq12}
\ed
where $\bW=(W_1,\cdots,W_n)^T$. Note that throughout this paper, whenever two vectors $\ba$ and $\bb$ are of the same length, $\ba \bb$ denotes the componentwise product. Given a predefined threshold value $\tau_n$, we select a set of variables as follows:
 \bg
  \mathcal{M}_{\tau_n}=\{1\leq j\leq p: \hat{u}_{nj}\geq\tau_n\}.
 \label{eq13}
 \ed

Alternatively, we can rank the covariates by the residual sum of squares of marginal nonparametric regressions, which is defined as
 \bg
 \hat{v}_{nj} = \| \bY - \hat{a}_{nj}(\bW) - \hat{b}_{nj}(\bW) \bX_j \|_n^2,
\label{eq14}
 \ed
 and we select variables as follows,
  \bg
  \mathcal{M}_{\nu_n}=\{1\leq j\leq p: \hat{v}_{nj}\leq \nu_n\},
\label{eq15}
 \ed
 where $\nu_n$ is a predefined threshold value.

It is worth noting that ranking by marginal utility $\hat{u}_{nj}$ is equivalent to ranking by the measure of goodness of fit  $\hat{v}_{nj}$. To see the equivalence, first note that
\bg
  \|\hat{a}_{nj}(\bW)+\hat{b}_{nj}(\bW)\bX_j\|^2_n=\frac1n\bY^T\bQ_{nj}(\bQ^T_{nj}\bQ_{nj})^{-1}\bQ^T_{nj}\bY,
\label{eq16}
 \ed
 and
  \bg
 \frac1n\sum\limits_{i=1}^nY_i(\hat{a}_{nj}(W_i)+\hat{b}_{nj}(W_i)X_{ji}) =\frac1n\bY^T\bQ_{nj}(\bQ^T_{nj}\bQ_{nj})^{-1}\bQ^T_{nj}\bY.
\label{eq17}
 \ed
It follows from (\ref{eq16}) and (\ref{eq17}) that
\bg
\hat{v}_{nj}=\|\bY\|^2_n- \|\hat{a}_{n0}(\bW)\|^2_n - \hat{u}_{nj}.
\label{eq18}
\ed
Since the first two terms on the right hand side of (\ref{eq18}) do not vary in $j$, ranking by $\hat{u}_{nj}$ is the same as that by $\hat{v}_{nj}$. Therefore, selecting variables with large marginal utility is the same as picking those that yield small marginal residual sum of squares.\\

To bridge $u_j$ and $\hat{u}_{nj}$, we define the population version of the marginal regression using B-spline basis. From now on, we will omit the argument in $\bB(W)$ and write $\bB$ whenever the context is clear.
Let $\tilde{a}_{j}(W)= \bB \tilde{\boldeta}_j$ and $\tilde{b}_{j}(W)= \bB \tilde{\btheta}_j$, where $\tilde{\boldeta}_{j}$ and $\tilde{\btheta}_{j}$ are the minimizer of
\bg
\min\limits_{\sbeta_{j},\sbtheta_{j}\in \mathbb{R}^{L_n}} \E[(Y-\bB \boldeta_{j}-\bB \btheta_{j} X_j)^2],
\label{eq19}
\ed
and $\tilde{a}_{0}(W) = \bB  \tilde{\boldeta}_0$, where $\tilde{\boldeta}_{0}$ is the minimizer of
\bg
\min\limits_{\sbeta_{0} \in \mathbb{R}^{L_n}} \E[(Y-\bB \boldeta_{0})^2].
\label{eq20}
\ed
It can be seen that
\bg
(\tilde{a}_j (W), \tilde{b}_j(W))^T &=& \diag(\bB, \bB)(\E[\bQ^T_j \bQ_j])^{-1}\E[\bQ^T_j Y], \label{eq21} \\
\tilde{a}_0(W) &=& \bB (\E[\bB^T\bB])^{-1} \E[\bB^T Y], \label{eq22}
\ed
where $\bQ_j = (\bB, X_j\bB)$ 
\bg
\tilde{u}_{j} &=& \|  \tilde{a}_j(W) + \tilde{b}_j(W)X_j \|^2 - \| \tilde{a}_0(W) \|^2 \label{eq23}  \nonumber \\
&=& \E[Y \bQ_j](\E[\bQ^T_j \bQ_j])^{-1}\E[\bQ^T_j Y] - \E[Y\bB](\E[\bB^T \bB])^{-1} \E[\bB^T Y].
\label{eq24}
\ed

\section{Sure Screening}
 In this section, we establish the sure screening properties of the proposed method for model (\ref{eq1}). Recall that by (\ref{eq6}) the population version of marginal utility quantifies the relationship between $X_j's$ and $Y$ as follows:
 \bg
 u_j = \E\left[ \frac{(\Cov[X_j,Y|W])^2}{\Var[X_j|W]}\right], \quad j=1,\cdots,p.
  \label{eq25}
 \ed
 Then the following two conditions guarantee that the marginal signal of the active components $\{u_j\}_{j\in \mathcal{M}_*}$ does not vanish.

\begin{description}
 \item [(i)] Suppose for $j=1,\cdots,p$, $\Var[X_j|W]$ is uniformly bounded away from 0 and infinity on $\mathcal{W}$, where $\mathcal{W}$ is the compact support of $W$. That is, there exist some positive constants $h_1$ and $h_2$, such that $0<h_1\leq \Var[X_j|W]\leq h_2<\infty$.

\item[(ii)] $\min_{j\in\mathcal{M}_*}\E[(\Cov[X_j,Y|W])^2] \geq c_1L_nn^{-2\kappa}$, for some $\kappa>0$ and $c_1>0$.
\end{description}
Then under conditions (i) and (ii),
 \bg
  \min_{j\in\mathcal{M}_*}u_j \geq c_1 L_nn^{-2\kappa}/h_2.
  \label{eq26}
 \ed
 Note that in condition (ii), the number of basis functions $L_n$ is not intrinsic. By the Remark \ref{remark1} below,  $L_n$ should be chosen in correspondence to the smoothness condition of the nonparametric component. Therefore, condition (ii) depends only on $\kappa$ and smoothness parameter $d$ in condition (iii). We keep $L_n$ here to make the relationship more explicit.

\subsection{Sure Screening Properties}
The following conditions (iii)-(vii) are required for the B-spline approximation in marginal regressions and establishing the sure screening properties.
\begin{description}
 \item[(iii)] The density function $g$ of $W$ is bounded away from zero and infinity on $\mathcal{W}$. That is, $0<T_1\leq g(W)\leq T_2<\infty$  for some constants $T_1$ and $T_2$.
 \item[(iv)] Functions $\{a_j\}_{j=0}^p$ and  $\{b_j\}_{j=1}^p$  belong to a class of functions $\mathcal{B}$, whose $r$th derivative $f^{(r)}$ exists and is Lipschitz of order $\alpha$. That is,
 \bgn
  \mathcal{B}=\{f(\cdot): |f^{(r)}(s)-f^{(r)}(t)|\leq M|s-t|^\alpha \mbox{ for } s,t\in\mathcal{W}\},
 \edn
 for some positive constant $M$, where $r$ is a nonnegative integer and $\alpha\in(0,1]$ such that $d=r+\alpha>0.5$.
\item[(v)] Suppose for all $j=1,\cdots,p$, there exists a positive constant $K_1$ and $r_1\geq 2$, such that
\bg
 \mbox{P}(|X_j|>t|W) \leq \exp(1-(t/K_1)^{r_1}),
\label{eq29}
\ed  uniformly on $\mathcal{W}$, for any $t \geq 0$. Furthermore, let $m(\bX^*)=\E[Y|\bX,W]$, where $\bX^*=(\bX^T,W)^T$. Suppose there exists some positive constants $K_2$ and $r_2$ satisfying $r_1r_2/(r_1+r_2)\geq1$, such that
 \bg
  \mbox{P}(|m(\bX^*)|>t | W)\leq \exp(1-(t/K_2)^{r_2}).
  \label{eq30}
 \ed
 uniformly on $\mathcal{W}$, for any $t \geq 0$.

\item[(vi)] The random errors $\{\varepsilon_i\}_{i=1}^n$ are i.i.d with conditional mean 0, and there exists some positive constants $K_3$ and $r_3$ satisfying $r_1r_3/(r_1+r_3)>1$, such that
 \bg
  \mbox{P}(|\varepsilon|>t|W)\leq \exp(1-(t/K_3)^{r_3}),
  \label{eq31}
 \ed
  uniformly on $\mathcal{W}$, for any $t \geq 0$.

\item [(vii)] There exists some constant $\xi\in(0,1/h_2)$ such that $L_n^{-2d-1}\leq c_1(1/h_2-\xi)n^{-2\kappa}/M_1$.
\end{description}

\begin{prop}\label{prop1}
 Under conditions (i)-(v), there exists a positive constant $M_1$ such that
\bg
u_{j} - \tilde{u}_{j} \leq M_1L_n^{-2d}.
 \label{eq27}
\ed
In addition, when $L_n^{-2d-1}\leq c_1(1/h_2-\xi)n^{-2\kappa}/M_1$ for some
$\xi\in(0,1/h_2)$, we have
\bg
 \min_{j\in\mathcal{M}_*}\tilde{u}_{j}\geq c_1\xi L_n n^{-2\kappa}.
 \label{eq28}
\ed

 \end{prop}

\begin{remark}\label{remark1}
It follows from Proposition \ref{prop1} that the minimum signal level of $\{\tilde{u}_{j}\}_{j\in\mathcal{M}_*}$ is approximately the same as $\{u_{j}\}_{j\in\mathcal{M}_*}$, provided that the approximation error is negligible. It also shows that the number of basis functions $L_n$ should be chosen as
$$L_n \geq Cn^{2\kappa/(2d+1)},$$ for some positive constant $C$. In other words, the smoother the underlying function is (i.e., the larger $d$ is), the smaller  $L_n$  we can take.
\end{remark}

The following Theorem \ref{thm1} provides the sure screening properties of the nonparametric independence screening method proposed in Section \ref{2.2}.

\begin{theorem} \label{thm1}
Suppose conditions (i)-(vi) hold.
\begin{enumerate}
\item [(i)]  If $n^{1-4\kappa} L_n^{-3}\to \infty$ as $n\to \infty$, then for any $c_2>0$, there exist some positive constants $c_3$ and $c_4$ such that
\bg
 &&P\left(\max_{1\leq j\leq p}|\hat{u}_{nj} - \tilde{u}_j|\geq c_2L_nn^{-2\kappa}\right) \nonumber\\
&\leq&12 p_n L_n \{(2+L_n)\exp(-c_3n^{1-4\kappa}L^{-3}_n)+3L_n\exp(-c_4L^{-3}_n n)\}.
 \label{eq32}
\ed

\item [(ii)] If condition (vii) also holds, then by taking $\tau_n=c_5L_nn^{-2\kappa}$ with $c_5=c_1\xi/2$, there exist positive constants $c_6$ and $c_7$ such that
\bg
P\left(\mathcal{M_*}\subset\hat{\mathcal{M}}_{\tau_n}\right)  & \geq & 1-12s_n L_n\{(2+L_n)\exp(-c_6n^{1-4\kappa}L^{-3}_n) \nonumber \\
&& +3L_n\exp(-c_7L^{-3}_n n)\}.
 \label{eq33}
\ed
\end{enumerate}
\end{theorem}

\begin{remark} \label{remark2}
According to Theorem \ref{thm1} ,  we can handle NP dimensionality
\bgn
p=o(\exp \{ n^{1-4 \kappa} L_n^{-3}\}).
\edn
It shows that the number of spline bases $L_n$ also affects the order of dimensionality: the smaller $L_n$ is, the higher dimensionality we can handle. On the other hand, Remark \ref{remark1} points out that it is required $L_n\geq Cn^{2\kappa/(2d+1)}$ to have a good bias property. This means that the smoother the underlying function is (i.e. the larger $d$ is), the smaller $L_n$ we can take, and consequently higher dimensionality can be handled.  The compatibility of these two requirements requires that $\kappa < (d+0.5)/(4d+5)$, which implies that $\kappa < 1/4$.  We can take $L_n=O(n^{1/(2d+1)})$, which is the optimal convergence rate for nonparametric regression (Stone, 1982).  In this case, the allowable  dimensionality can be as high as
\bgn
p=o(\exp \{ n^{\frac{2(d-1)}{2d+1}} \} ).
\edn
\end{remark}

 \subsection{False Selection Rates}
 According to (\ref{eq28}), the ideal case for vanishing false-positive rate is when
 \bgn
  \max\limits_{j\notin\mathcal{M}_*} \tilde{u}_j =o(L_nn^{-2\kappa})
 \edn
so that there is a natural separation between important and unimportant variables. By Theorem \ref{thm1}(i), when (\ref{eq32}) tends to zero, we have with probability tending to 1 that
 \bgn
  \max\limits_{j\notin\mathcal{M}_*} \hat{u}_{nj} \leq c L_nn^{-2\kappa},\mbox{ for any } c>0.
 \edn
Consequently, by choosing $\tau_n$ as in Theorem \ref{thm1}(ii), NIS can achieve the model selection consistency under this ideal situation, i.e.,
\bgn
 P\left(\hat{\mathcal{M}}_{\tau_n}=\mathcal{M}_*\right)=1-o(1).
\edn
In particular, this ideal situation occurs under the partial orthogonality condition, i.e., $\{X_j\}_{j\in\mathcal{M}_*}$ is independent of $\{X_i\}_{i \notin \mathcal{M}_*}$ given $W$, which implies $u_j = 0$ for $j \not \in\mathcal{M}_*$

In general, the model selection consistency can not be achieved by a single step of marginal screening.  The marginal probes can not separate important variables from unimportant variables. The following Theorem \ref{thm2} quantifies how the size of selected models is related to the matrix of basis functions and the thresholding parameter $\tau_n$.

\begin{theorem}\label{thm2}
Under the same conditions in Theorem \ref{thm1}, for any $\tau_n=c_5L_nn^{-2\kappa}$, there exist positive constants $c_8$ and $c_9$ such that
\bg
P\left\{|\hat{\mathcal{M}}_{\tau_n}|\leq O(n^{2\kappa}\lambda_{\max}(\bSigma))\right\} &\geq& 1-12 p_n L_n \Big \{(2+L_n)\exp(-c_8 n^{1-4\kappa} L^{-3}_n )\nonumber \\
 && + 3  L_n \exp(-c_9nL_n^{-3}) \Big \},
\label{eq34}
\ed
where $\bSigma=\E [\bQ^T\bQ]$, and $\bQ=(\bQ_1,\cdots,\bQ_{p})$ is a functional vector of $2p_nL_n$ dimension.
\end{theorem}

\section{Iterative Nonparametric Independence Screening}
As \cite{Fan2008} points out, in practice the nonparametric
independence screening (NIS) would still suffer from false negative
(i.e., miss some important predictors that are marginally weakly
correlated but jointly correlated with the response), and false
positive (i.e., select some unimportant predictors which are highly
correlated with the important ones). Therefore, we adopt an
iterative framework to enhance the performance of this method. We
repeatedly apply the large-scale variable screening (NIS) followed
by a moderate-scale variable selection, where we use group-SCAD
penalty as our selection strategy. In the NIS step, we propose two
methods to determine a data-driven threshold for screening, which
result in Conditional-INIS and Greedy-INIS, respectively.

\subsection{Conditional-INIS Method}
The conditional-INIS method builds upon \textit{conditional} random
permutation in determining the thresholding $\tau_n$. Recall the
random permutation used in \cite{Fan2011}, which generalizes that
\cite{Zhao2010}.  Randomly permute $\bY$ to get $\bY_{\bpi}=
(Y_{\pi_1}, \cdots, Y_{\pi_n})^T$ and compute $\hat{u}_{nj}^{\bpi}$,
where $\bpi$ is a permutation of $\{1, \cdots, n\}$, based on the
randomly coupled data $\{({Y}_{\pi_i}, W_i, \bX_i)\}_{i=1}^n$ that
has no relationship between covariates and response.  Thus, these
estimates serve as the baseline of the marginal utilities under the
null model (no relationship).  To control the false selection rate
at $q/p$ under the null model, one would choose the screening
threshold be $\tau_q$, the $q$th-ranked magnitude of
$\{\hat{u}_{nj}^{\bpi},~j=1,\cdots,p\}$.  Thus, the NIS step selects
variables $\{j:\hat{u}_{nj} \geq \tau_q\}$.  In practice, one
frequently uses $q = 1$, namely, the largest marginal utility under
the null model.

When the correlations among covariates are large, there will be hardly any differentiability between the marginal utilities of the true variables and the false ones. This makes the selected variable set very large to begin with and  hard to proceed the rest of iterations with limited false positives. For numerical illustrations, see section \ref{conditional section}. Therefore, we propose a \textit{conditional} permutation method to tackle this problem. Combining the other steps, our Conditional-INIS algorithm proceeds as follows.

\begin{description}\label{inis1}
\item[0.] For $j=1,\cdots,p$, compute $$\hat{u}_{nj} = \|\hat{a}_{nj}(\bW) + \hat{b}_{nj}(\bW) \bX_j \|_n^2 - \|\hat{a}_{n0}(\bW) \|_n^2,$$ where
the estimates are defined in (\ref{eq9}) and (\ref{eq10}) using $\{(\bY, \bW, \bX_j), j=1, \cdots, p \}$.
Select the top $K$ variables by ranking their marginal utilities $\hat{u}_{nj}$, resulting in the index subset $\mathcal{M}_0$ to condition upon.

\item[1.]
Regress $\bY$ on $\{(\bW,\bX_j), j\in \mathcal{M}_0\}$, and get intercept $\hat{\beta}_{n0}(W) $ and their functional coefficients' estimators $\{\hat{\beta}_{nj}(W), j\in \mathcal{M}_0\}$.
Conditioning on $\mathcal{M}_0$, the $n$-dimensional partial residual is
\bgn
 \bY^*=\bY-\hat{\beta}_{n0}(\bW) - \sum_{ j\in \mathcal{M}_0} \bX_j \hat{\beta}_{nj}(\bW).
\edn
For all $j\in\mathcal{M}^c_0$, compute $\hat{u}_{nj}^*$
using $\{(\bY^*, \bW, \bX_j), j \in \mathcal{M}^c_0\}$,
which measures the additional utility of each covariate conditioning on the selected set $\mathcal{M}_0$.

To determine the threshold for NIS, we apply random permutation on the partial residual $\bY^*$, which yields $\bY^*_{\bpi}$. Compute $\hat{u}^{* \bpi}_{nj}$ based on the decoupled data $\{(\bY^*_{\bpi}, \bW, \bX_j), j \in \mathcal{M}^c_0\}$.
Let $\tau^*_q$ be the $q$th-ranked magnitude of $\{\hat{u}^{* \bpi}_{nj}, j\in\mathcal{M}^c_0\}$.  Then, the active variable set of variables is chosen as
\bgn
 \mathcal{A}_1=\{j:\hat{u}^*_{nj}  \geq \tau^*_q, j\in\mathcal{M}^c_0\}\cup\mathcal{M}_0.
\edn
In our numerical studies,  $q = 1$.

\item[2.] Apply the  group-SCAD penalty on $\mathcal{A}_1$ to select a subset of variables $\mathcal{M}_1$.
Details about the implementation of SCAD will be described later.

\item[3.] Repeat step 1-2, where we replace $\mathcal{M}_0$ in step 1 by $\mathcal{M}_l$, $l=1,2,\cdots$, and get  $\mathcal{A}_{l+1}$ and $\mathcal{M}_{l+1}$ in step 2. Iterate until $\mathcal{M}_{l+1}=\mathcal{M}_k$ for some $k\leq l$ or $|\mathcal{M}_{l+1}|\geq \zeta_n$, for some prescribed positive integer $\zeta_n$.

\end{description}

\subsection{Greedy-INIS Method}
Following \cite{Fan2011}, we also implemented a greedy version of
INIS method. We skip step 0 and start from step 1 in the algorithm
above (i.e., take $\mathcal{M}_0=\emptyset$), and select the top
$p_0$ variables that have the largest marginal norms $\hat{u}_{nj}$.
This NIS step is followed by the same group-SCAD penalized
regression as in step 2. We then iterate these steps until there are
two identical subsets or the number of variables selected exceeds a
prespecified $\zeta_n$. In our simulation studies,  $p_0$ is set as
$1$.


\subsection{Implementation of SCAD}
In the group-SCAD step, variables are selected as $\mathcal{M}_{l} =
\{j \in \mathcal{A}_l: \hat\bgamma_j^{(l)} \neq \textbf{0} \} $
through minimizing the following objective function: \bg
\min\limits_{\sbgamma_0,\sbgamma_j \in \mathbb{R}^{L_n} }
\frac{1}{n} \sum_{i=1}^n \Bigl (Y_i - \bB(W_i)\bgamma_0 - \sum_{j\in
\mathcal{A}_l} \bB (W_i) X_{ji} \bgamma_j \Bigr )^2 + \sum_{j\in
\mathcal{A}_l} p_{\lambda} (||\bgamma_j||_B), \label{eq35} \ed where
$||\bgamma_j||_B = \sqrt{\frac{1}{n} \sum_{i=1}^n (\sum_{k=1}^{L_n}
B_{jk}(W_i)\gamma_{jk})^2},$ and $p_{\lambda}(\cdot)$ is the SCAD
penalty such that \bgn p'_{\lambda}(|x|)&=&\lambda I (|x|\leq
\lambda) + \frac{(a \lambda - |x|)_+}{a-1} I (|x|>\lambda), \edn
with $p_{\lambda}(0)=0$. We set $a=3.7$ as suggested and solve the
optimization above via local quadratic approximations
\citep{Fan2001}. $\lambda$ is chosen by BIC criteria $n\log
(\hat{\sigma}^2_\epsilon) + k L_n \log n$, where $k$ is the number
of covariates chosen. By \cite{Antoniadis2001} and \cite{Yuan2006},
the norm-penalty in (\ref{eq35}) encourages the group selection.

\section{Numerical Studies}
In this section, we carry out several simulation studies to assess the performance of our proposed methods. If not otherwise stated, the common setup for the following simulations are: cubic B-spline, $L_n=7,$ sample size $n=400$ , the number of variables $p=1000$, and the number of simulations  $N=200$ for each example.

\subsection{Comparison of Minimum Model Size}
In this study, as in \cite{Fan2010}, we illustrate the performance
of NIS method in terms of the minimum model size (MMS) needed to
include all the true variables, i.e., to possess sure screening
property.

\paragraph{\textit{Example 1}}
Following \cite{Fan2010}, we first consider a linear model as a
special case of the varying coefficient model. Let
$\{X_k\}_{k=1}^{950}$  be i.i.d. standard normal random variables
and
\[
X_k=\sum_{j=1}^s (-1)^{j+1}X_j/5 + \sqrt{1-\frac{s}{25}} \xi_k, \quad k=951, \cdots, 1000,
\]
where $\{\xi_k\}_{k=951}^{1000}$ are standard normal random variables. We construct the following model: $ Y= \bbeta ^{T} \bX + \epsilon $, where $\epsilon \sim \mathcal{N} (0, \sqrt{3}^2)$ and $ \bbeta= (1,-1,1,-1,\cdots)^{T} $ has $s$ nonzero components. To carry out NIS, we define an exposure $W$ independently from the standard uniform distribution.

\begin{figure}[h]
\centering
\includegraphics[width=\textwidth]{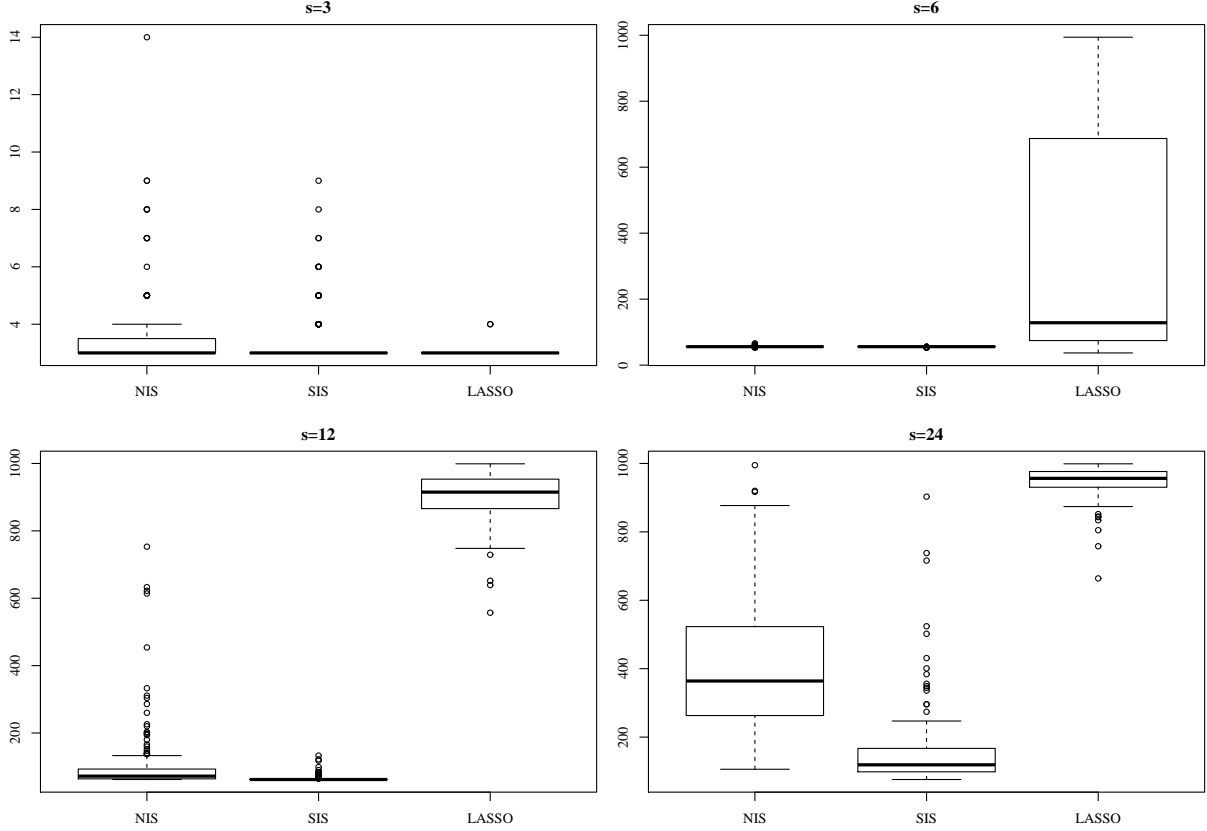}
\caption{\label{ex1} Boxplots of minimum model sizes (left to right: NIS, Lasso and SIS) for Example 1 under different true models. }
\end{figure}

We compare NIS, Lasso and SIS (independence screening for linear models). The boxplots of minimum model size are presented in Figure \ref{ex1}. Note that when $s>5$, the irrepresentable condition fails, and Lasso performs badly even in terms of pure screening. On the other hand, SIS performs better than NIS because the coefficients are indeed constant, and there are fewer parameters ($p$) involved in SIS than those of NIS ($pL_n$).

\paragraph{\textit{Example 2}}
For the second example, we illustrate that when the underlying model's coefficients are indeed varying, we do need nonparametric independence screening. Let $ \{U_1, U_2, \cdots U_{p+2}\} $ be i.i.d. uniform random variables on $ [0,1] $, based on which we construct $\bX$ and $W$  as follows:
\bgn
X_j &=& \frac{ U_j + t_1 U_{p+1}}{1+t_1}, \quad j=1,\cdots,p, \qquad
W =  \frac{ U_{p+2} + t_2 U_{p+1}}{1+t_2},
\edn
where $t_1$ and $t_2$ controls the correlation among the covariates $\bX$ and the correlation between $\bX$ and $W$, respectively. When $t_1=0$, $X_j$'s are uncorrelated, and when $t_1=1$ the correlation is $0.5$. If $t_1=t_2=1$, $X_j$'s and $W$ are also correlated with correlation coefficient 0.5.

For the varying coefficients part, we take coefficient functions
\bgn
\beta_1 (W) = W, \quad \beta_2 (W) = (2W-1)^2,\quad \beta_3 (W) = \sin(2 \pi W).
\edn
The true data generation model is
\[
Y = 5 \beta_1 (W) \cdot X_1 + 3 \beta_2 (W) \cdot X_2 + 4 \beta_3 (W) \cdot X_3 + \epsilon,
\]
where $\epsilon$'s are i.i.d. standard Gaussian random variable. \\

\begin{figure}[h]
\centering
\includegraphics[width=\textwidth]{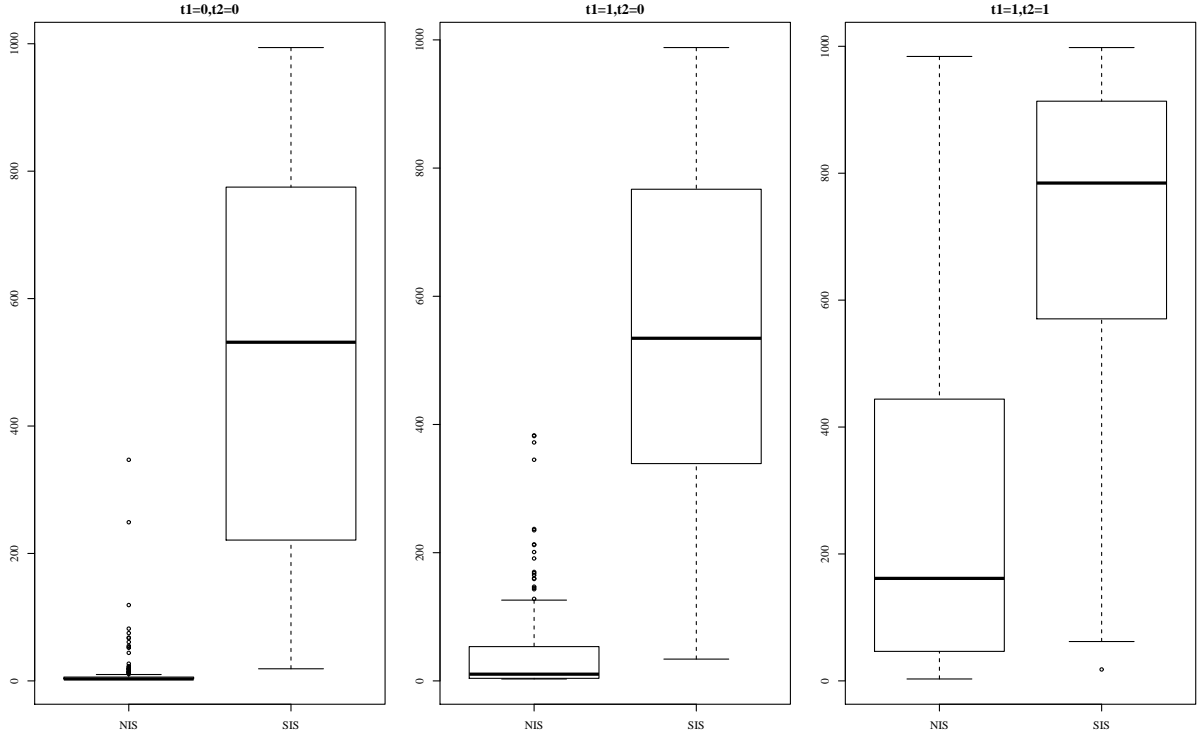}
\caption{\label{ex2} Boxplots of minimum model sizes (left: NIS, right: SIS) for Example 2 under different correlation settings.}
\end{figure}

Under different correlation settings, the comparison MMS between NIS and SIS methods are presented in Figure \ref{ex2}. When the correlation gets stronger, independence screening becomes harder.

\subsection{Comparison of Permutation and Conditional Permutation}{\label{conditional section}}
In this section, we illustrate the performance the conditional random permutation method.

\paragraph{\textit{Example 3}}
Let $\{Z_1,\cdots,Z_p\}$ be i.i.d. standard normal, $\{U_1,U_2\}$ be i.i.d. standard uniformly distributed random variables, and the noise $\epsilon$ follows the standard normal distribution. We construct $\{W,\bX\} $ and $Y$ as follows:
\bgn
X_j&=&\frac {Z_j+t_1 U_1} {1+t_1}, j=1,\cdots, p, \qquad
W =  \frac{U_2 + t_2 U_1}{1+t_2},\\
Y&=& 2X_1 + 3W \cdot X_2+(W+1)^2 \cdot X_3 +\frac{4\sin(2\pi W)}{2-\sin(2\pi W)} \cdot X_4+ \epsilon.
\edn
We will take $t_1 = t_2 = 0$, resulting in uncorrelated case and $t_1=3$ and $t_2=1$, corresponding to $\mbox{corr}(X_j, X_k ) = 0.43$ for all $j \not = k$ and
$\mbox{corr}(X_j, W) = 0.46$.
By taking $q=1$ (i.e., take the maximum value of the marginal utility of the permuted estimates), we report the average of the true positive number (TP), model size, the lower bound of the marginal signal of true variables and the upper bound of the marginal signal of false variables for different correlation settings based on $200$ simulations.  Their robust standard deviations are also reported therein.



Based on Table \ref{taa2}, we see that when the correlation gets stronger, although sure screening properties can be achieved most of the time via
unconditional ($K=0$) random permutation thresholding, the model size becomes very large and therefore the false selection rate is high. The reason is that there is no differentiability between the marginal signals of the true variables and the false ones. This drawback makes the original random permutation not a feasible method to determine the screening threshold in practice.

\begin{table}\caption{\label{taa2}Model size and marginal signals under different correlation settings (Example 3)}
\centerline{
\begin{tabular}{cc|c|c|c|c|c}
\hline\hline
\multicolumn{2}{c|}{\multirow{2}{*}{Model}}  & \multirow{2}{*}{TP} & \multirow{2}{*}{ Size}  & \multirow{2}{*}{$\min\limits_{j \in \mathcal{M}^* \backslash \mathcal{M}_0} \hat{u}_{nj}^*$} &\multirow{2}{*}{$\max\limits_{j \in \mathcal{M}^{*c}\backslash \mathcal{M}_0} \hat{u}_{nj}^* $} & \multirow{2}{*}{$\max\limits_{j \in \{1,\cdots, p\}\backslash \mathcal{M}_0} \hat{u}_{nj}^{*\bpi}$}\\
  &&&&&&\\
  \hline \multirow{2}{*}{K=0} &  $ t_1=0, t_2=0 $      & $4.00(0)$ & $6.68(2.99)$  & $2.96(0.72)$ & $1.22(0.18)$  & 1.12(0.15)\\
    &   $ t_1=3, t_2=1$    & $4.00(0)$ & $886.49(88.81)$  & $0.61(0.10)$ & $0.58(0.07)$  & 0.22(0.03)\\
    \hline \multirow{2}{*}{K=1} &  $ t_1=0, t_2=0 $      & $4.00(0)$ & $5.70(1.49)$  & $2.83(0.57)$ & $0.75(0.10)$  & 0.72(0.11)\\
    &   $ t_1=3, t_2=1$    & $4.00(0)$ & $202.50(154.85)$  & $0.28(0.06)$ & $0.20(0.03)$ & 0.11(0.02)\\
\hline \multirow{2}{*}{K=4} &  $ t_1=0, t_2=0 $      & $4.00(0)$ & $5.14(1.49)$  & NA & $0.06(0.01)$ & 0.06(0.01)\\
    &   $ t_1=3, t_2=1$    & $4.00(0)$ & $4.98(0.75)$  & $0.16(0.05)$ & $0.05(0.01)$ & 0.06(0.01)\\
  \hline \multirow{2}{*}{K=8} &    $ t_1=0, t_2=0 $      & $4.00(0)$ & $8.92(0.75)$  & NA & $0.05(0.01)$  & 0.05(0.01)\\
    &   $ t_1=3, t_2=1$    & $3.99(0)$ & $8.43(0.75)$  & $0.11(0.03)$ & $0.04(0.01)$  & 0.05(0.01)\\
        \hline\hline
\end{tabular}
}
\end{table}

We now applied the conditional permutation method, whose performance is illustrated in Table \ref{taa2} for a few choices of tuning parameter $K$.
The screening threshold is taken as $\tau_q$ with $q=1$.
Generally speaking, although the lower bound of the true positives' signals may be smaller than the upper bound of false variables' signals, the largest $K$ norms still have a high possibility to contain at least some true variables. When conditioning on this small set of more relevant variables, the marginal contributions of false positives get weaker. Note that in the absence of correlation, when $K \geq s$ (here $s=4$), the first $K$ variables have already included all the true variables (i.e., $\mathcal{M}^* \backslash \mathcal{M}_0 = \emptyset $), hence the minimum of true signal is not available. In other cases, we see that the gap between the  marginal signals of true variables and false variables become large enough to differentiate them. Table \ref{taa2} shows that by using the thresholding via the conditional permutation method, not only the sure screening properties are still maintained, but also the model sizes are dramatically reduced.

\subsection{Comparison of Model Selection and Estimation}
In this section we explore the performance of Conditional-INIS and Greedy-INIS method. In our iterative framework, conditional permutation serves as the initialization step (step 0) and we take $K=5$ in the rest of the paper.
For each method, we report the average number of true positive (TP), false positive (FP), prediction error (PE), and their robust standard deviations. Here the prediction error is the mean squared error calculated on the test dataset of size $n/2=200$ generated from the same model. As a measure of the complexity of the model, signal-to-noise-ratio (SNR), defined by ${\var(\bbeta^T(W) \bX)}/{\var(\epsilon)}$, is computed.
Table~\ref{tab3} reports the results using the simulated model specified in Example 3.  We now illustrate the performance by using another example.

\begin{table}\caption{\label{tab3} Average values of the number of true positives (TP), false positives (FP), and prediction error (PE) for simulated model in Example 3. Robust standard deviations are given in parentheses. \newline}
\centering
        \begin{tabular}{c|c|c|c|c|c|c|c|c}
        \hline  \hline
   {Model} &\multicolumn{2}{|c|}{Correlation} & \multicolumn {3}{|c|}{Conditional-INIS} &  \multicolumn {3}{|c}{Greedy-INIS} \\
      \cline{2-9} &  X's & X's-W  &    TP& FP&PE &TP &FP&PE \\
    \hline  $ t_1=0, t_2=0$  &$0$&$0$ & $4$ &  $0.54$&$1.10$& $4$&$13.01$&$1.41$\\
              $(\SNR \approx 16.85)$&&&$(0)$&$(0.75)$&$(0.05)$&$(0)$&$(3.73)$&$(0.17)$\\

    \hline  $ t_1=2, t_2=0$  &$0.25$&$0$ & $4$ &  $0.20$&$0.78$& $ 4$&$0.41 $&$1.10$\\
                   $(\SNR \approx 3.66)$&&&$(0)$&$(0)$&$(0.06)$&$(0)$&$(0)$&$(0.05)$\\

      \hline  $ t_1=2, t_2=1$  & $0.25$ & $0.36$ & $3.97$ & $0.26$&$1.27$& $3.90$&$0.14 $&$1.63$\\
                          $(\SNR \approx 3.21)$&&&$(0)$&$(0)$&$(0.24)$  &$(0)$&$(0)$&$(0.41)$\\

    \hline  $ t_1=3, t_2=0$  & $0.43$ & $0$ & $4$ &  $0.19$&$1.03$& $ 3.99$&$ 0.57$&$1.22$\\
              $(\SNR \approx 3.32)$&&&$(0)$&$(0)$&$(0.06)$      &$(0)$&$(0)$&$(0.07)$\\

  \hline    $ t_1=3, t_2=1$ &$0.43$&$0.46$    & $3.95$  & $0.31$& $1.30$ &$ 3.77 $ &$ 0.27 $&$1.29$ \\
                 $(\SNR \approx 2.81)$&&&$(0)$&$(0.75)$&$(0.12)$ &$(0)$&$(0)$&$(0.17)$\\
            \hline\hline
        \end{tabular}
\end{table}

\paragraph{\textit{Example 4}}
Let $\{W,\bX\} $ , $Y$ and $\epsilon$ be the same as in \textit{Example 3}. We now introduce more complexities in the following model:
\bgn
Y &=& 3W \cdot X_1 + (W+1)^2 \cdot X_2 + (W-2)^3\cdot X_3 + 3(\sin(2 \pi W))\cdot X_4 \\
   && + \exp(W)\cdot X_5 +2 \cdot X_6 + 2 \cdot X_7+ 3\sqrt{W}  \cdot X_8 + \epsilon.
\edn
The results are present in Table \ref{tab4}.

\begin{table}\caption{\label{tab4} Average values of the number of true positives (TP), false positives (FP), and prediction error (PE) for the model in Example 4. Robust standard deviations are given in parentheses. \newline}

\centering
        \begin{tabular}{c|c|c|c|c|c|c|c|c}
        \hline  \hline
   {Model} &\multicolumn{2}{|c|}{Correlation} & \multicolumn {3}{|c|}{Conditional-INIS} &  \multicolumn {3}{|c}{Greedy-INIS} \\
      \cline{2-9} &  X's & X's-W  &    TP& FP&PE &TP &FP&PE \\
    \hline  $ t_1=0, t_2=0$  &$0$&$0$ & $ 8$ & $ 0.21$ &$1.24$   &$8$ & $10.71$ & $1.57$ \\
             $(\SNR \approx 47.68)$&&&$ (0)$ & $ (0)$ &$(0.09) $&$ (0) $ & $ (3.73)$ & $ (0.20) $\\

    \hline  $ t_1=2, t_2=0$  &$0.25$&$0$ & $ 8$ & $ 0.13 $ & $1.17 $& $8$&$0.60$ &$1.16$\\
                 $(\SNR \approx 9.40)$&&&$ (0)$ & $ (0)$ & $(0.09)$& $ (0)$&$ (0)$ &$ (0.10)$\\

      \hline  $ t_1=2, t_2=1$  & $0.25$ & $0.36$ &  $ 7.80$ & $0.20$ & $2.16$&$ 7.55 $ &$ 0.26 $ &$ 2.26$ \\
                            $(\SNR \approx 8.62)$&&&$ (0)$ & $ (0)$ & $(0.58)$& $ (0.75)$&$ (0)$ &$(0.70) $\\

    \hline  $ t_1=3, t_2=0$  & $0.43$ & $0$ &$ 7.90$ & $ 0.10$ &  $ 1.21$&$ 7.98 $&$ 0.71  $ &$1.29$\\
                     $(\SNR \approx 8.18)$&&&$(0) $ & $ (0)$ & $(0.12)$& $ (0)$&$ (0)$ &$ (0.10)$\\

  \hline    $ t_1=3, t_2=1$ &$0.43$&$0.46$    &  $7.75$ & $0.18$ &  $1.65$  & $ 7.35 $&$0.28$ &$1.84$  \\
                       $(\SNR \approx 7.61)$&&& $(0)$ & $(0)$ &  $(0.26)$  & $(0.75)$&$(0)$ &$(0.42)$\\
            \hline\hline
        \end{tabular}
\end{table}

Through the examples above, Conditional-INIS and Greedy-INIS show comparable performance in terms of TP, FP and PE. When the covariates are independent or weakly correlated, sure screening is easier to achieve and false positive is rare; as the correlation gets stronger, we see a decrease in TP and an increase in FP. It seems that  Greedy-INIS selects slightly more false positives than Conditional-INIS, the reason being that in each step Greedy-INIS selects the top variable(s) by fitting the residuals conditional on previously chosen variable set and tends to overfit. However, the coefficient estimates for these false positives are fairly small, hence they do not affect prediction error very much. Regarding computation efficiency, Conditional-INIS performs better in our simulated examples, as it usually only requires two to three iterations, while Greedy-INIS would  need at least $s/p_0$ iterations (here $p_0=1$ and $s=4$ and $8$ respectively for Examples 3 and 4).

\subsection{Real Data Analysis on Boston Housing Data}
In this section we illustrate the performance of our method through
a real data analysis on Boston Housing Data \citep{Harrison1978}. 
This dataset contains housing data for 506 census tracts of Boston
from the 1970 census. Most empirical results for the housing value
equation are based on a common specification \citep{Harrison1978},
\bgn
\log (\mbox{MV}) &=& \beta_0 + \beta_1 \mbox{ RM}^2 + \beta_2 \mbox{AGE} + \beta_3 \log(\mbox{DIS}) + \beta_4 \log(\mbox{RAD}) + \beta_5 \mbox{TAX} \\
          &&+ \beta_6 \mbox{PTRATIO} + \beta_7 (\mbox{B}-0.63)^2 + \beta_8 \log(\mbox{LSTAT}) + \beta_{9} \mbox{CRIM} \\
          &&+ \beta_{10} \mbox{ZN} + \beta_{11} \mbox{INDUS} + \beta_{12} \mbox{CHAS} + \beta_{13} \mbox{NOX}^2 + \epsilon,
\edn
where the dependent variable MV is the median value of owner-occupied homes, the independent variables are quantified measurement of its neighborhood whose description can be found in the manual of R package \textit{mlbench}. The common specification uses $\mbox{RM}^2$ and $\mbox{NOX}^2$ to get a better fit, and for comparison we take these transformed variables as our input variables.

\begin{figure}[h]
\centering
\includegraphics[width=\textwidth]{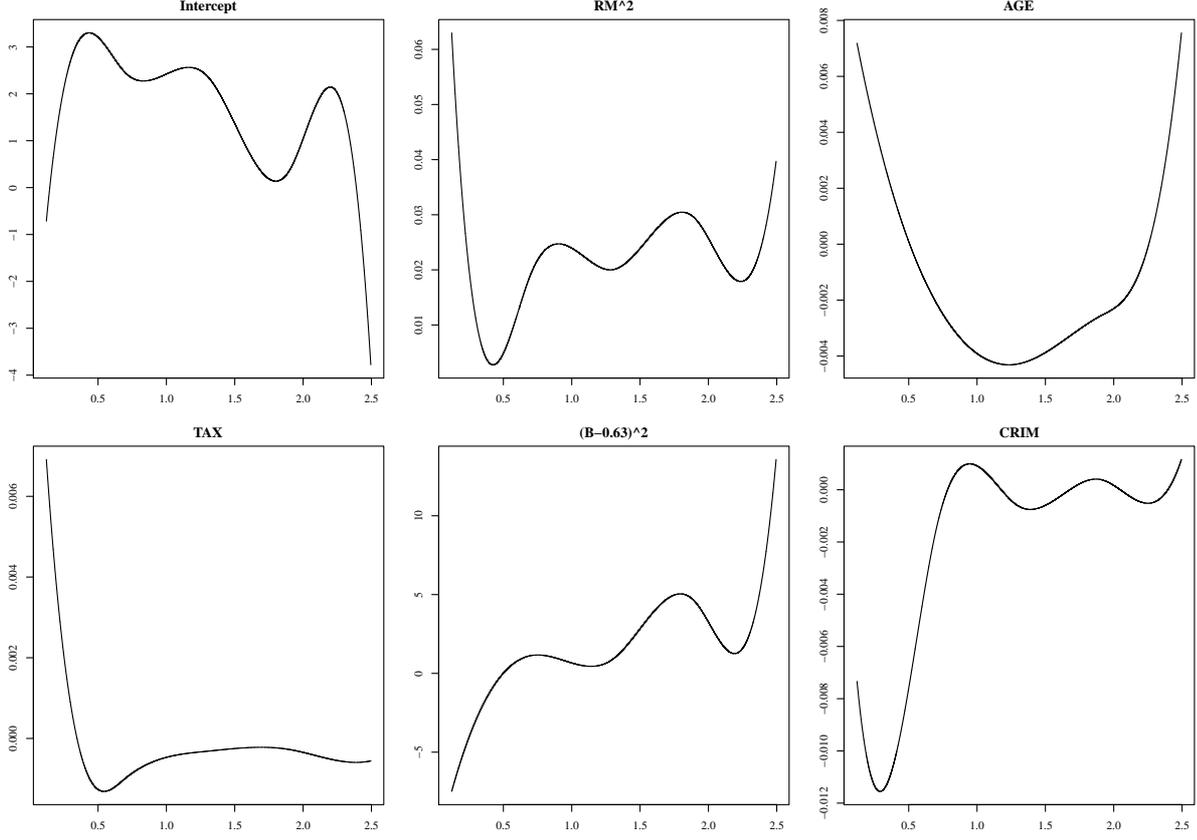}
\caption{\label{c_scad} Fitted functional estimates $\hat{\beta}_j(W) 's$ selected by Conditional-INIS.}
\end{figure}

To exploit the power of varying coefficient model, we take the variable $\log(\mbox{DIS})$, the weighted distances to five employment centers in the Boston region, as the exposure variable. This allows us to examine how the distance to the business hubs interact with other variables. It is reasonable to assume that the impact of other variables on housing price varies with the distance, which is an important characteristic of the neighborhood, i.e. the geographical accessibility to employment. Interestingly, Conditional-INIS selects the following submodel:
\bg
\log (\mbox{MV}) &=& \beta_0(W) + \beta_1 (W) \cdot \mbox{RM}^2  + \beta_2 (W) \cdot \mbox{AGE}+ \beta_5(W) \cdot \mbox{TAX}  \nonumber\\
&&+ \beta_7(W) \cdot (\mbox{B}-0.63)^2 + \beta_{9}(W) \cdot \mbox{CRIM} + \epsilon,
\label{realdata_CINIS}
\ed
where $W = \log(\mbox{DIS})$. The estimated functions $\hat{\beta}_j (W)$'s are presented in Figure \ref{c_scad}. This varying coefficient model shows very interesting aspects of housing valuation.  The evidence of nonlinear interactions with the accessibility is clearly evidenced. For example, RM is the average number of rooms in owner units, which represents the size of a house. Therefore, the marginal cost of a big house is higher in employment centers where population is concentrated and supply of mansions is limited.  The cost per room decreases as one moved away from the business centers and then gradually increases. CRIM is the crime rate in each township, which usually has a negative impact, and from its varying coefficient we see that it is a bigger concern near (demographically more complex) business centers.  AGE is the proportion of owner units built prior to 1940, and its varying coefficient has a parabola shape: positive impact on housing values near employment centers and suburb areas, while negative effects in between. NOX (air pollution level) is generally a negative impact, and the impact is larger when the house is near employment centers where air is presumably more polluted than suburb area.

We now evaluate the performance of our INIS method in a high
dimensional setting. To accomplish this, let $\{Z_1,\cdots,Z_p\}$ be
i.i.d. the standard normal random variables and $U$ follow the
standard uniform distribution. We then expand the data set by adding
the artificial predictors: \bgn X_j = \frac {Z_j+ t U} {1+t},
j=s+1,\cdots, p. \edn Note that $\{W, X_1,\cdots, X_s\}$ are the
independent variables in original data set ($s=13$ here) and the
variables $\{X_j\}_{j=s+1}^p$ are known to be irrelevant to the
housing price, though the maximum spurious correlation of these 987
artificial predictors to the housing price is now small. We take
$p=1000$, $t=2$, and randomly select $n=406$ samples as training
set, and compute prediction  mean squared  error (PE) on the rest
$100$ samples. As a benchmark for comparison, we also do regression
fit on $\{W, X_1,\cdots, X_s\}$ directly using SCAD penalty without
screening procedure. We repeat $N = 100$ times and report the
average prediction error and model size, and their robust standard
deviation. Since $\{X_j\}_{j=s+1}^p$ are artificial variables, we
also include the number of artificial variables selected by each
method as a proxy for false positives. The results are presented in
Table \ref{tab5}.

\begin{table}\caption{\label{tab5} Prediction error (PE) , model size and selected noise variables (SNV) over 100 repetitions and their robust standard deviations (in parentheses) for
Conditional-INIS (p = 1000), Greedy-INIS (p = 1000), and SCAD fit (p=12).}

\centering
        \begin{tabular}{c|c|c|c}
         \hline
  method & PE & Size & SNV \\
      Conditional-INIS (p = 1000)   &$ 0.046 (0.020)$ & $5.55(0.75) $ & $0 (0) $ \\

       Greedy-INIS (p = 1000)  & $ 0.048 (0.020)$&$ 4.80 (1.49)$ &$ 0.01 (0)$\\

         SCAD fit (p=12)   & $ 0.052 (0.019)$&$ 6.05 (1.87)$ &$ NA $\\
     \hline
        \end{tabular}
\end{table}

As seen from Table \ref{tab5}, our methods are very effective in filtering noise variables in a high dimensional setting, and can achieve comparable prediction error as if the noise were absent. In conclusion, the proposed INIS methodology is very useful in high-dimensional scientific discoveries, which can select a parsimonious close-to-truth model and reveal interesting relationship between variables, as illustrated in this section.

\section*{Acknowledgements}
This project was supported by
the National Institute of General Medical Sciences of the National Institutes of Health through Grant Numbers R01-GM072611 and R01-GMR01GM100474 and National Science Foundation grant DMS-1206464. The bulk of the research was carried while Yunbei Ma was a postdoctoral fellow at Princeton University.

\section*{Appendix}
\subsection*{A.1. Properties of  B-splines}
Our estimation use the B-spline basis, which have the following properties (de Boor 1978):
For each $j=1,\cdots, p$ and $k=1,\cdots,L_n$, $B_{k}(W)\geq 0$ and $\sum_{k=1}^{L_n}B_{k}(W)=1$ for $W\in\mathcal{W}$. In addition, there exist positive constants $T_3$ and $T_4$ such that for any $\eta_{k}\in\mathbb{R},~k=1,\cdots,L_n$,
\bg
 L^{-1}_nT_3\sum_{k=1}^{L_n}\eta^2_{k}\leq \int\left(\sum_{k=1}^{L_n}\eta_{k} B_{k}(w)\right)^2dw\leq L^{-1}_nT_4\sum_{k=1}^{L_n}\eta^2_{k}.
 \label{bspline1}
\ed
Then under condition (iii), there exist positive constants $C_1$ and $C_2$ such that for $k=1,\cdots,L_n$,
\bg
  C_1L_n^{-1}\leq \E[ B^2_{k}(W)]\leq C_2L_n^{-1},
  \label{bspline2}
\ed
where $C_1=T_1T_3$ and $C_2=T_2T_4$.

Furthermore, under condition (iii), it follows from (\ref{bspline1}) that for any $\boldeta=(\eta_1,\cdots,\eta_{L_n})^T\in\mathbb{R }^{L_n}$ such that $\|\boldeta\|_2^2=1$,
\bgn
  C_1L^{-1}_n\leq\boldeta^T \E[\bB^T\bB]\boldeta\leq C_2 L^{-1}_n.
\edn
Or equivalently,
\bg
C_1L^{-1}_n\leq\lambda_{min}(\E[\bB^T\bB])\leq\lambda_{max}(\E[\bB^T\bB])\leq C_2L^{-1}_n.
\label{bspline3}
\ed

\subsection*{A.2. Technical Lemmas}
Some technical lemmas needed for our main results are shown as follows. Lemma \ref{lemma1} and Lemma \ref{lemma2} give some characterization of exponential tails , which becomes handy in our proof. Lemma \ref{lemma3} and Lemma \ref{lemma4} is a Bernstein type inequality.

\begin{lemma}\label{lemma1}
Let $X$, $W$ be random variables. Suppose $X$ has a conditional exponential tail: $\mbox{P}(|X|>t|W)\leq \exp (1-(t/K)^{r})$ for all $t\geq0$ and uniformly on the compact support of $W$, where $K > 0 $ and $r\geq1$. Then for all $m\geq2$,
\bg
\E(|X|^m|W)\leq emK^mm!.
\label{le01}
\ed
\end{lemma}
{\bf Proof}. Recall that for any non-negative random variable $Z$, $\E[Z|W]=\int_0^\infty \mbox{P}\{Z\geq t|W\}dt$. Then we have
\bgn
 \E(|X|^m|W)&=&\int_0^\infty \mbox{P}\{|X|^m\geq t|W\}dt\\
 &\leq&\int_0^\infty \exp(1- (t^{1/m}/K)^r)dt\\
 &= & \frac{emK^m}{r}\Gamma(\frac{m}{r}).
\label{eq36}
\edn
The lemma follows from the fact $r \geq 1$.

\begin{lemma}\label{lemma2}
Let $Z_1$, $Z_2$ and $W$ be random variables. Suppose that
there exist $K_1$, $K_2>0$ and $r_1$, $r_2\geq1$ such that $r_1r_2/(r_1+r_2)\geq1$, and
$$\mbox{P}(|Z_i|>t|W)\leq \exp (1-(t/K_i)^{r_i}),\quad i=1,2$$
for all $t\geq0$ and uniformly on $\mathcal{W}$.
Then for some $r^*\geq 1$ and $K^*>0$,
\bg
 \mbox{P}(|Z_1Z_2|>t|W)\leq \exp (1-(t/K^*)^{r^*})
 \label{eq37}
\ed
for all $t\geq0$ and uniformly on $\mathcal{W}$.
\end{lemma}
{\bf Proof}. For any $t>0$, let $M=(tK_2^{r_2/r_1}/K_1)^{\frac{r_1}{r_1+r_2}}$ and $r=r_1r_2/(r_1+r_2)$. Then uniformly on $\mathcal{W}$, we have
\bgn
 \mbox{P}(|Z_1Z_2|>t|W)&\leq&\mbox{P}(M|Z_1|>t|W)+\mbox{P}(|Z_2|>M|W)\\
  &\leq&\exp\{1-(t/K_1M)^{r_1}\}+\exp\{1-(M/K_2)^{r_2}\}\\
  &=&2\exp\{1-(t/K_1K_2)^{r}\}.
\edn
Let $r^*\in[1,r]$ and $K^*=\max\{(r^*/r)^{1/r}K_1K_2,~(1+\log2)^{1/r}K_1K_2\}$. It can be shown that $G(t)=(t/K_1K_2)^r-(t/K^*)^{r^*}$ is increasing when $t>K^*$. Hence $G(t)> G(K^*) \geq \log2$ when $t>K^*$, which implies when $t>K^*$,
\bgn
 \mbox{P}(|Z_1Z_2|>t|W)\leq 2\exp\{1-(t/K_1K_2)^{r_1}\}\leq \exp\{1-(t/K^*)^{r^*}\}.
\edn
On the other hand, when $t\leq K^*$,
\bgn
 \mbox{P}(|Z_1Z_2|>t|W)\leq 1\leq \exp\{1-(t/K^*)^{r^*}\}.
\edn
Lemma \ref{lemma2} holds.

\begin{lemma}\label{lemma3}
(Bernstein inequality, lemma 2.2.11, \cite{van1996}). For
independent random variables $Y_1,\cdots,Y_n$ with mean zero such
that $\E [|Y_i|^m] \leq m!M^{m-2}\nu_i/2$ for every $m\geq2$ (and
all $i$) and some constants $M$ and $\nu_i$. Then \bgn
 P(|Y_1+\cdots+Y_n|>x)\leq 2\exp\{-x^2/(2(\nu+Mx))\},
\edn
for $v\geq\nu_1+\cdots+\nu_n$.
\end{lemma}

\begin{lemma}\label{lemma4}
(Bernstein's inequality, lemma 2.2.9, \cite{van1996}). For
independent random variables $Y_1,\cdots,Y_n$ with bounded range
$[-M,M]$ and mean zero, \bgn
 P(|Y_1+\cdots+Y_n|>x)\leq 2\exp\{-x^2/(2(\nu+Mx/3))\},
\edn for $\nu\geq\var(Y_1+\cdots+Y_n)$.
\end{lemma}

The following lemmas are needed for the proof of Theorem \ref{thm1}.

\begin{lemma}\label{lemma5}
Suppose conditions (i) and (iii)-(vi) hold. For any $\delta>0$, there exist some positive constants $b_1$ and $b_2$ such that for $j=1,\cdots,p$, $k=1,\cdots,L_n$,
\bgn
 P\left(\left|\frac1n\sum_{i=1}^nX_{ji}B_{k}(W_i)Y_i-\E [X_j B_{k}Y] \right| \geq \frac{\delta } {n} \right) \leq 4 \exp\left \{-\frac{\delta^2}{b_1L_n^{-1}n + b_2 \delta)}\right \},
\edn
 and
  \bgn
P\left(\left|\frac1n\sum_{i=1}^nB_{k}(W_i)Y_i-\E [B_{k}Y] \right| \geq \frac{\delta}{  n} \right) \leq 4 \exp\left \{-\frac{\delta^2}{b_1L_n^{-1}n + b_2 \delta }\right \}.
\edn
\end{lemma}

{\bf Proof}.  Recall $m(\bX_i^*) = E(Y_i|\bX_i, W_i)$.  Let $Z_{jki}=X_{ji}B_{k}(W_i)m(\bX^*_i) - \E [X_jB_{k}(W)m(\bX^*)]$ and $\xi_{jki}=X_{ji}B_{k}(W_i)\varepsilon_i$. Then
\bgn
&&\left |\frac1n\sum_{i=1}^n X_{ji}B_{k}(W_i)Y_i-\E [X_jB_{k}(W)Y] \right | \\
&=& \left |\frac1n\sum_{i=1}^n \Bigl ( X_{ji}B_{k}(W_i)m(\bX^*_i) - \E [X_jB_{k}(W)m(\bX^*)] +X_{ji}B_{k}(W_i)\varepsilon_i \Bigr ) \right | \\
&\leq & \left |\frac1n \sum_{i=1}^n Z_{jki}\right | + \left |\frac1n \sum_{i=1}^n \xi_{jki} \right |.
\edn

We first bound $\frac1n \sum_{i=1}^n Z_{jki}$. Note that for each $j$ and
$k$,   $\{Z_{jki}\}_{i=1}^n$ are a sequence of independent random variables with mean zero. By condition (v), (\ref{bspline2}),  and Lemmas \ref{lemma1} and \ref{lemma2}, we have for every $m \geq 2$,  there exists a constant $K_4>0$, such that
\bg
\E|Z_{jki}|^m &\leq& 2^m \E |X_{ji}B_{k}(W_i) m(\bX^*_i)|^m  \nonumber\\
&\leq& 2^m  \E[B^m_{k}(W_i) \E[|X_{ji}m(\bX^*_i)|^m|W_i]] \nonumber \\
&\leq & 2^m \E[B^2_{jk}(W_i) emK^m_4m!] \nonumber \\
& \leq & m! (2K_4)^{m-2} (8emK_4^2C_2L_n^{-1})/2, \label{jf}
\ed
where the first inequality comes from the Minkowski inequality. Hence, it follows from Lemma \ref{lemma3} that for any $\delta >0$,
\bg
P\left(\Bigl |\frac1n \sum_{i=1}^n Z_{jki}\Bigr | \geq \frac{\delta}{2n} \right) \leq 2 \exp \left \{ -\frac{\delta^2 }{ 64 em K_4^2C_2 L_n^{-1}n + 8 K_4\delta}\right \}
\label{eq38}
\ed

Next we bound $\frac1n \sum_{i=1}^n \xi_i$. Again $\xi_i$'s are centered independent random variables. By conditions (v)-(vi), (\ref{bspline2}), and Lemmas \ref{lemma1} and \ref{lemma2}, we have for every $m \geq 2$, there exists a constant $K_5>0$, such that
\bgn
\E |\xi_i|^m &=& \E [B_{k}^m(W_i) \E[|X_{ji}\varepsilon_i|^m|W_i]] \\
& \leq & m! K_5^{m-2} (2emK_5^2C_2L_n^{-1})/2.
\edn
Thus, according to Lemma \ref{lemma3},
\bg
P\left(\Big |\frac1n \sum_{i=1}^n \xi_i \Big | \geq \frac{\delta}{ 2n} \right) \leq 2 \exp
 \left \{ -\frac{\delta^2 }{16emK_5^2C_2L_n^{-1}n + 4K_5\delta} \right \}.
\label{eq39}
\ed

Similarly, we can show that
 \bg
&&P\left( \Big |\frac1n \sum_{i=1}^n B_{k}(W_i)m(\bX^*_i)-\E [B_{k}(W)m(\bX^*)]
   \Big | \geq \frac{\delta}{2n} \right) \nonumber\\
&\leq& 2 \exp \left \{ -\frac{\delta^2}{64 em K_2^2C_2 L_n^{-1}n + 8 K_2\delta} \right \}
\label{eq40}
\ed
and
 \bg
P\left(\Big |\frac1n \sum_{i=1}^n B_{k}(W_i)\varepsilon_i \Big | \geq \frac{\delta}{ 2n} \right)
\leq 2 \exp \left \{ -\frac{\delta^2}{16emK_3^2C_2L_n^{-1}n +4K_3\delta} \right \}.
\label{eq41}
\ed

Let $b_1 = 16 e m C_2 \max (4 K_4^2, K_5^2, 4 K_2^2, K_3^2)$ and $b_2 = \max (8K_4, 4K_5, 8K_2, 4K_3)$. Then, the combination of (\ref{eq38}) - (\ref{eq41}) by union bound of probability yields the desired result.  $\Box$

\begin{lemma}\label{lemma6} Under conditions (i), (iii) and (v), there exist positive constants $C_3$ and $C_4$, such that for $j = 1, \cdots, p$,
\bg
C_3L^{-1}_n\leq\lambda_{min}(\E[\bQ^T_j\bQ_j])\leq\lambda_{max}(\E[\bQ^T_j\bQ_j])\leq C_4L^{-1}_n.
\label{bspline4}
\ed
\end{lemma}

{\bf Proof}. Recall that $\bQ_{j}=(\bB, X_j \bB)$. For any $\boldeta = (\boldeta_1^T,\boldeta^T_2)^T
\in\mathbb{R }^{2L_n}$ such that $\| \boldeta \|_2^2=1$,
\bgn
\boldeta^T \E[\bQ^T_j\bQ_j] \boldeta & = & \E \left[(\bB\boldeta_1,\bB\boldeta_2) \left(\begin{array}{cc} 1 & \E[X_j|W] \\ \E[X_j|W] & \E[X^2_j|W] \end{array}\right) \left(\begin{array}{c} \bB\boldeta_1 \\ \bB\boldeta_2\end{array}\right)\right].
\edn
Consider eigenvalues $\lambda_1$ and $\lambda_2$ ($\lambda_1 > \lambda_2$) of the $2 \times 2$ middle matrix on the right hand side of the equation above, we have $\lambda_1 + \lambda_2 = 1 + \E[X_j^2|W]$ (trace) and $\lambda_1 \cdot \lambda_2 = \Var[X_j|W]$ (determinant). Therefore, by Lemma~\ref{lemma1}
$$
  \lambda_1 \leq 1 + \E[X_j^2|W] \leq 1 + 4 e  K_1^2
$$
and by assumption (i)
$$
    \lambda_2 \geq \frac{\Var[X_j|W]}{\E[X_j^2|W]+1} \geq \frac{h_1}{1 + 4 e  K_1^2}.
$$
Using the above two bounds on the minimum and maximum eigenvalues, we have
\bgn
\frac{h_1}{1+4 e K_1^2} \E[(\bB\boldeta_1)^2 + (\bB\boldeta_2)^2 ] \leq \boldeta^T \E[\bQ^T_j\bQ_j] \boldeta \leq ( 1+4 e K_1^2) \E[(\bB\boldeta_1)^2 + (\bB\boldeta_2 )^2].
\edn
By (\ref{bspline3}), we have
\bgn
\frac{h_1C_1}{1+4 e K_1^2}L^{-1}_n \leq \boldeta^T \E[\bQ^T_j\bQ_j] \boldeta \leq ( 1+4 e K_1^2)C_2L^{-1}_n.
\edn
Take $C_3 = h_1C_1L^{-1}_n/(1+4 e K_1^2)$ and $C_4 = (1+4 e K_1^2)C_2L^{-1}_n$, result follows. \\

Throughout the rest of the proof, for any matrix $\bA$, let $\|\bA\|=\sqrt{\lambda_{\max}(\bA^T\bA)}$ be the operator norm and $\|\bA\|_\infty=\max_{i,j}|A_{ij}|$ be the infinity norm.

\begin{lemma}\label{lemma7}
Suppose conditions (i), (iii) and (v) hold. For any $\delta>0$ and $j=1,\cdots,p$, there exist some positive constants $b_3$ and $b_4$ such that
\bgn
P\left(\Big \|\frac1n\bQ^T_{nj}\bQ_{nj}-\E [\bQ^{T}_j\bQ_j]\Big \|\geq L_n\delta/n\right)
\leq 6L^2_n\exp\left \{- \frac{\delta^2}{b_3L_n^{-1}n + b_4 \delta} \right \},
\edn
and
\bgn
P \left(\Big \| \frac1n\bB^T_{n}\bB_{n}-\E [\bB^{T}\bB]\Big  \| \geq L_n\delta/n\right) \leq 6L^2_n\exp \left \{-
   \frac{\delta^2}{b_3L_n^{-1}n + b_4 \delta} \right \}.
\edn
 In addition, for any given positive constant $b_5$, there exists some positive constant $b_6$ such that
 \bgn
P\left(\left| \Big \|(\frac1n\bQ^T_{nj}\bQ_{nj})^{-1}\|-\|(\E [\bQ^{T}_j \bQ_j])^{-1}\Big \|\right|\geq b_5\|(\E [\bQ^{T}_j \bQ_j])^{-1}\|\right) \leq6L^2_n\exp\{-b_6L^{-3}_n n\},
 \edn
and for any positive constant $b_7$, there exists some positive constant $b_8$ such that
\bgn
P\left(\left | \Big \|(\frac1n\bB^T_{n}\bB_{n})^{-1} \Big \|-
   \Big \|(\E [\bB^{T} \bB])^{-1} \Big \|\right|\geq b_7\|(\E [\bB^{T} \bB])^{-1}\|\right) \leq6L^2_n\exp\{-b_8L^{-3}_n n\}.
 \edn
\end{lemma}

{\bf Proof}. Observe that for $j = 1, \cdots, p,$
\bgn
  \frac1n\bQ^T_{nj} \bQ_{nj}- \E[\bQ^{T}_j\bQ_j] =\left(\begin{array}{cc}{\bD}_{1}&\bD_{2j}\\
\bD^T_{2j}&\bD_{3j}\end{array}\right),
\edn
where
$\bD_{1}=\frac1n \sum\limits_{i=1}^{n} \bB^T(W_i)\bB(W_i)- \E[\bB^T\bB]$,
$\bD_{2j}=\frac1n \sum\limits_{i=1}^{n} X_{ji}\bB^T(W_i)\bB(W_i) - \E[X_j\bB^T\bB] $ and
$\bD_{3j}=\frac1n \sum\limits_{i=1}^{n} X^2_{ji}\bB^T(W_i)\bB(W_i) - \E[X^2_j\bB^T\bB].$
Then
 \bg
  \| \frac1n\bQ^T_{nj} \bQ_{nj}- \E[\bQ^{T}_j\bQ_j] \| &\leq&
  2L_n\| \frac1n\bQ^T_{nj} \bQ_{nj}- \E[\bQ^{T}_j\bQ_j]  \|_\infty\nonumber\\
  &=& 2L_n\max ( \|\bD_{1}\|_\infty,~\|\bD_{2j}\|_\infty,~\|\bD_{3j}\|_\infty).
 \label{eq42}
 \ed

We first bound $\|\bD_{1}\|_\infty$. Recall that $0\leq B_k(\cdot) \leq 1$ on $\mathcal{W}$, so
 \bgn
| B_k(W_i) B_l(W_i) - \E[B_k(W)B_l(W)] |
\leq 2,
 \edn
for all $k$ and $l$ 
By (\ref{bspline2}),
 \bgn
  \Var \left(B_k(W_i) B_l(W_i)  - \E[B_k(W)B_l(W)] \right) \leq \E [B_k^2(W) B_l^2(W)] \leq C_2L^{-1}_n.
 \edn
By Lemma \ref{lemma4}, we have
 \bgn
 & & P\left(|\frac1n\sum\limits_{i=1}^n B_k(W_i) B_l(W_i)  - \E[B_k(W)B_l(W)] |
  \geq  \delta/6n\right) \\
 & \leq & 2\exp\{-\delta^2/(72C_2L^{-1}_nn + 24 \delta)\}.
 \edn
It then follows from the union bound of probability that
 \bg
  P\left(\|\bD_{1}\|_\infty\geq \delta/6n\right)\leq2L^2_n\exp\{-\delta^2/(72 C_2L^{-1}_nn+24\delta)\}.
  \label{eq43}
 \ed

We next bound $\|\bD_{2j}\|_\infty$. Note that for $k,l = 1, \cdots, L_n$,
\bgn
&& \E [ |X_{ji}  B_k(W_i) B_l(W_i) - \E[X_j B_k(W) B_l(W)] |^m ]\\
 & \leq &  2^m \E[ |X_{ji}B_{k}(W_i)B_{l}(W_i)|^m ]\\
 &\leq & 2^m  \E [ |X_{ji} B_{k}(W_i)|^m ] \\
&=& 2^m  \E[ \E[|X_{ji}|^m|W_i] B_{k}^m(W_i) ]\\
&\leq& m! (2K_1)^{m-2} (8 emK_1^2C_2L_n^{-1})/2,
\edn
where Lemma~\ref{lemma1} was used in the last inequality.
By Lemma \ref{lemma3}, we have
\bgn
&& P\left(| \frac1n\sum\limits_{i=1}^nX_{ji}  B_k(W_i) B_l(W_i) - \E[X_j B_k(W) B_l(W)] |\geq \delta /6n \right) \\
&\leq& 2 \exp \{ -\delta^2 / (576 emK_1^2C_2L_n^{-1}n + 24 K_1\delta)\}.
\edn
 It then follows from the union bound of probability
that
 \bg
  P\left(\|\bD_{2j}\|_\infty\geq \delta/6n\right)\leq2L^2_n\exp\{-\delta^2/ (576 emK_1^2C_2L_n^{-1}n + 24 K_1\delta)\}.
  \label{eq44}
 \ed

Similarly we can bound $\|\bD_{3j}\|_\infty$. For every $m \geq 2$, for $k,l = 1, \cdots,L_n$, there exists a constant $K_6>0$ such that
\bgn
&& \E[ |X_{ji}^2B_{k}(W_i)B_{l}(W_i) - \E [ X_j^2B_{k}(W)B_{l}(W)] |^m ]\\
&\leq & 2^m \E[ \E[|X_{ji}^2|^m|W_i] B_{k}^m(W_i)]\\
&\leq& m! (2K_6)^{m-2} (8emK_6^2C_2L_n^{-1})/2. \edn
By Lemma \ref{lemma3}, we have
\bgn
&& P\left(|X_{ji}^2B_{k}(W_i)B_{l}(W_i) - \E [ X_j^2B_{k}(W)B_{l}(W)]|\geq \delta /6n \right) \\
&\leq& 2 \exp \{ -\delta^2 / (576 emK_6^2C_2L_n^{-1}n+ 24  K_6 \delta)\}.
\edn
It then follows from the union bound of probability that
\bg
  P\left(\|\bD_{3j}\|_\infty\geq \delta/6n\right)\leq2L^2_n\exp\{-\delta^2/(576 emK_6^2C_2L_n^{-1}n + 24 K_6 \delta)\}.
  \label{eq45}
\ed

Let $b_3=72 C_2\max\{1 , 8emK^2_1,8em K^2_6\}$ and $b_4=24 \max\{1, K_1, K_6\}$, then combining (\ref{eq42})-(\ref{eq45}) we have
 \bg
  P\left(\Big \| \frac1n\bQ^T_{nj} \bQ_{nj}- \E[\bQ^{T}_j\bQ_j] \Big \|\geq L_n\delta/n\right)\leq 6L^2_n\exp\left \{-\frac{\delta^2}{b_3L_n^{-1}n + b_4 \delta} \right \}.
  \label{eq46}
 \ed
Observe that $\| \frac1n\bB^T_{n} \bB_{n}- \E[\bB^{T}\bB] \| \leq 2L_n \|\bD_{1}\|_\infty$. Thus,  we have also proved that
 \bg
  P\left(\Big \| \frac1n\bB^T_{n} \bB_{n}- \E[\bB^{T}\bB] \Big \|\geq L_n\delta/n\right)\leq 6L^2_n\exp\left \{- \frac{\delta^2}{b_3L_n^{-1}n + b_4 \delta}
  \right \}.
  \label{eq47}
 \ed

We next prove the second part of the lemma. Note that for any
symmetric matrices $\bA$ and $\bB$ \citep{Fan2011},
 \bg
 |\lambda_{\min}(\bA)-\lambda_{\min}(\bB)|\leq \max\{|\lambda_{\min}(\bA-\bB)|,|\lambda_{\min}(\bB-\bA)|\}.
 \label{eq48}
 \ed
It then follows from (\ref{eq48}) that
 \bgn
\left |\lambda_{\min}(\frac1n\bQ^T_{nj}\bQ_{nj})-\lambda_{\min}(\E[\bQ^{T}_j\bQ_j])
\right |
&\leq& 2 L_n\Big \|\frac1n\bQ^T_{nj}\bQ_{nj} - \E [\bQ^{T}_j\bQ_j]\Big \|_{\infty},
\edn
which implies that
 \bg
&& P\left(\left |\lambda_{\min}(\frac1n\bQ^T_{nj}\bQ_{nj})-\lambda_{\min}(\E[\bQ^{T}_j\bQ_j])
\right |\geq L_n \delta /n \right)\nonumber\\
&\leq&6L^2_n\exp\{-\delta^2/(b_3L_n^{-1}n + b_4\delta)\}.
\label{eq49}
\ed
Let $\delta = b_9 C_3 L^{-2}_n n$ in (\ref{eq49}) for $b_9 \in (0,1)$. According to (\ref{bspline4}),
we have
\bg
 && P\left( \left |\lambda_{\min}(\frac1n\bQ^T_{nj}\bQ_{nj})-\lambda_{\min}( \E [\bQ^{T}_j\bQ_j] ) \right | \geq b_9 \lambda_{\min}( \E [\bQ^{T}_j\bQ_j])\right)\nonumber\\
 &\leq& 6L^2_n\exp(-b_6L^{-3}_nn),
 \label{eq50}
\ed
for some positive constant $b_6$. Next observe the fact that for $x,y>0, a \in (0,1) $ and $b = 1/(1-a) -1$,
\bgn
|x^{-1} - y^{-1}| \geq b y^{-1}  \mbox{ implies } |x-y| \geq a y.
\edn
This is because $x^{-1} - y^{-1} \geq b y^{-1}$ is equivalent to $x^{-1} \geq \frac{1}{1-a} y^{-1}$, or $x - y \leq -ay$; on the other hand, $x^{-1} - y^{-1} \leq b y^{-1}$ implies $x^{-1} \leq (1 - \frac{a}{1-a})y^{-1} \leq (1- \frac{a}{1+a}) y^{-1}$ as $a \in (0,1)$, and therefore $x - y \geq ay$.
Then let $b_5=1/(1-b_9)-1$, it follows from (\ref{eq50}) that
\bg
 &&P\left( \left |(\lambda_{\min}(\frac1n\bQ^T_{nj}\bQ_{nj}))^{-1}
  -(\lambda_{\min}(\E[\bQ^{T}_j\bQ_j)])^{-1}\right |
 \geq b_5(\lambda_{\min}(\E[\bQ^{T}_j\bQ_j]))^{-1}\right)\nonumber\\
&\leq& 6L^2_n\exp(-b_6L^{-3}_n n).
 \label{eq51}
\ed

Following the same proof, by (\ref{bspline3}) we also have for any positive constant $b_7$, there exists some positive constant $b_8$, such that
\bg
 &&P\left(\left |(\lambda_{\min}(\frac1n\bB^T_{n}\bB_{n}))^{-1}
 -(\lambda_{\min}(\E[\bB^{T}\bB)])^{-1} \right |
 \geq b_7(\lambda_{\min}(\E[\bB^{T}\bB]))^{-1}\right)\nonumber\\
&\leq& 6L^2_n\exp(-b_8L^{-3}_n n).
 \label{eq52}
\ed
The second part of the lemma then follows from the fact that for any symmetric matrix $\bA$, $\lambda_{\min}(\bA)^{-1}=\lambda_{\max}(\bA^{-1})$. $\Box$

\subsection*{A.3. Proof of Main Results}
{\bf Proof of Proposition \ref{prop1}.}
Note that $\E[Y|W,X_j] = a_j(W) + b_j(W)X_j$.
By \cite{Stone1982}, there exist $\{a^*_j \}_{j=0}^p$ and $\{b^*_j\}_{j=1}^p \in \mathcal{S}_n$ such that $\|a_j-a^*_j\|_\infty\leq M_2L^{-d}_n$ and $\|b_j-b^*_j\|_\infty\leq M_2L^{-d}_n$, where $\mathcal{S}_n$ is the space of polynomial splines of degree $l \geq 1$ with normalized B-spline basis $\{B_{k}, k=1,\cdots,L_n\}$, and $M_2$ is some positive constant. Here $\|\cdot\|_\infty$ denotes the sup norm. Let $\boldeta^*_j$ and $\btheta^*_j$ be $L_n$-dimensional vectors such that for $a^*_j (W)= \bB(W) \boldeta^*_j$ and $b^*_j (W)=\bB_j(W)\btheta^*_j$. \\
Recall that $\tilde{a}_{j}(W) =\bB(W)\tilde{\boldeta}_j$ and $\tilde{b}_{j}(W) = \bB(W)\tilde{\btheta}_j$.
By definition of $\tilde{\boldeta}_j$ and $\tilde{\btheta}_j$, we have
\bgn
(\tilde{a}_{j}, \tilde{b}_{j}) &=& \arg \min_{ a_{j}, b_{j} \in \mathcal{S}_n} \E[(Y - a_{j}(W) - b_{j}(W)X_j)^2] \\
&=& \arg \min_{a_{j}, b_{j}  \in \mathcal{S}_n} \E[(\E[Y|W,X_j] - a_{j}(W) - b_{j}(W)X_j)^2],
\edn
and therefore $\| \E[Y|W,X_j] - \tilde{a}_{j} - \tilde{b}_{j} X_j\|^2 \leq \| \E[Y|W,X_j] -a_{j}^* - b_{j}^*X_j\|^2$. In other words,
\bgn
\|\tilde{a}_{j}  + \tilde{b}_{j} X_j- (a_j + b_jX_j) \|^2 &\leq& \| (a_{j}^* + b_{j}^*X_j) -(a_j + b_jX_j) \|^2\\
&\leq& 2 \|a_j - a^*_j\|^2 + 2 \|(b_j-b^*_j)X_j\|^2 \\
& \leq & 2 M_2^2 L_n^{-2d} (1+ \E[X_j^2]).
\edn
On the other hand, by the least-squares property,
\bgn
\E[(Y - \tilde{a}_{j}  - \tilde{b}_{j} X_j)(\tilde{a}_{j}  +\tilde{b}_{j} X_j)] &=& 0,
\edn
and by conditioning in $W_j$ and $X_j$, we have
\bgn
\E[(Y-a_{j} - b_{j}X_j)(\tilde{a}_{j}  + \tilde{b}_{j} X_j)] &=& 0.
\edn
The last two equalities imply that
$$
\E[({a}_{j}  + {b}_{j} X_j - \tilde{a}_{j}  - \tilde{b}_{j} X_j)(\tilde{a}_{j}  +\tilde{b}_{j} X_j)] =0
$$
Thus, by the Pythagorean theorem, we have
\bgn
\|a_j + b_jX_j \|^2 = \|\tilde{a}_{j} + \tilde{b}_{j} X_j\|^2 +  \|\tilde{a}_{j} + \tilde{b}_{j} X_j - a_j- b_jX_j \|^2,
\edn
and
\bg
\|a_{j}+ b_{j}X_j\|^2 - \|\tilde{a}_{j} + \tilde{b}_{j} X_j \|^2  \leq 2 M_2^2 L_n^{-2d} (1+ \E[X_j^2]).
\label{eq53}
\ed
Similary, we have
\bg
\|a_{0} \|^2 -\| \tilde{a}_{0}\|^2 \leq  M_2^2 L_n^{-2d}.
\label{eq54}
\ed
By taking $M_1 = M_2^2 (8e K^2+3)$ (c.f. Lemma \ref{lemma1}), the first part of Proposition \ref{prop1} follows from (\ref{eq53}) and (\ref{eq54}):
\bg
u_j - \tilde{u}_{j} &=& \|a_{j}+ b_{j}X_j\|^2 - \|a_{0} \|^2  -( \|\tilde{a}_{j}  +\tilde{b}_{j} X_j \|^2 -\|\tilde{a}_{0}\|^2) \nonumber \\
&\leq& M_1 L_n^{-2d}.
\label{eq55}
\ed
By (\ref{eq26}) and (\ref{eq55}), we have
\bgn
\min_{j\in\mathcal{M}_*} \tilde{u}_{j} \geq c_1L_nn^{-2\kappa}/h_2-M_1L_n^{-2d}.
\edn
Then the desired result follows from $L_n^{-2d-1}\leq c_1(1/h_2-\xi)n^{-2\kappa}/M_1$ for some $\xi\in(0,1/h_2)$. $\Box$

{\bf Proof of Theorem \ref{thm1}}. We first prove part (1). Note that
\bgn
\hat{u}_{nj} - \tilde{u}_j &=& S_1 + S_2,
\edn
where $$S_1 = \frac 1n \bY^T \bQ_{nj} (\bQ^T_{nj} \bQ_{nj})^{-1} \bQ^T_{nj}\bY - \E[Y \bQ_j] (\E[\bQ^T_j \bQ_j])^{-1} \E[\bQ^T_jY],$$
and
$$S_2 = \frac 1n \bY^T \bB_{n} (\bB^T_{n} \bB_{n})^{-1} \bB^T_{n}\bY - \E[Y \bB] (\E[\bB^T\bB])^{-1} \E[\bB^T Y].$$

We first focus on $S_1$. Let $\ba_n =\frac 1n \bQ^T_{nj} \bY $,  $\ba=\E[\bQ^T_jY]$, $\bU_n = (\frac 1n \bQ^T_{nj} \bQ_{nj} )^{-1}$ and $\bU = (\E[\bQ^T_j \bQ_j])^{-1}$.  Then
\bgn
S_1 &=& \ba^T_n\bU_n\ba_n-\ba^T\bU\ba \\
&=&(\ba_n-\ba)^T{\bU}_n(\ba_n-\ba)+2(\ba_n-\ba)^T{\bU}_n\ba+\ba^T({\bU}_n-{\bU})\ba.
\edn
Denote the last three terms respectively by $S_{11}$, $S_{12}$, and  $S_{13}$.

We first deal with $S_{11}$. Note that
\bg
|S_{11}| \leq \|\bU_n\|\cdot\|\ba_n-\ba\|_2^2.
 \label{eq56}
\ed
 By Lemma \ref{lemma5} and the union bound of probability,
\bg
 P(\|\ba_n-\ba\|_2^2\geq 2L_n\delta^2 / n^2)
&\leq& 8L_n \exp\{-\delta^2/(b_1L_n^{-1}n+b_2\delta)\}.
 \label{eq57}
\ed
According to the second part of Lemma \ref{lemma7}, for any given positive constant $b_5$, there exists a positive constant $b_6$ such that
\bgn
  P\left(|\|{\bU}_n\|-\|{\bU}\||\geq b_5\|{\bU}\|\right)\leq  6L^2_n\exp\{ -b_6L^{-3}_nn \}.
\edn
Then it follows from (\ref{bspline4}) that
\bg
 P\left(\|{\bU}_n\|\geq (b_5+1)C^{-1}_3L_n\right)\leq 6L^2_n\exp \{ -b_6L^{-3}_nn \}.
 \label{eq58}
\ed
Combining (\ref{eq56})-(\ref{eq58}) and based on the union bound of probability, we have
\bg
&&P(|S_{11}|\geq2(b_5+1)C_3^{-1}L^2_n\delta^2/n^2)\nonumber\\
&\leq&8L_n\exp\{-\delta^2/(b_1L^{-1}_nn+b_2\delta)\}+6L^2_n\exp\{-b_6L^{-3}_nn\}.
\label{eq59}
\ed
We next bound $S_{12}$. Note that
\bg
 |S_{12}| \leq 2\|\ba_n-\ba\|_2\cdot\|{\bU}_n\|\cdot\|\ba\|_2
\label{eq60}
\ed
By Lemma~\ref{lemma1},
\bg
 \|\ba\|_2^2 &=&\| \E[\bB^T Y]\|_2^2+ \| \E[X_j\bB^T Y]\|_2^2 \nonumber\\
 &=&\sum_{k=1}^{L_n}( \E[B_{k}m(\bX^*)])^2+\sum_{k=1}^{L_n}( \E[X_jB_{k}m(\bX^*)])^2\nonumber\\
 &\leq& \sum_{k=1}^{L_n} ( \E[B_{k}^2m^2(\bX^*)]+ \E[B_{k}^2X_j^2m^2(\bX^*)]) \nonumber \\
 &\leq& 4eC_2(K_2^2+K^2_4),
 \label{eq61}
\ed
where the calculation as in (\ref{jf}) was used.

It follows from (\ref{eq57}), (\ref{eq58}),  (\ref{eq60}), (\ref{eq61}) and the union bound of probability that
\bg
 &&P(|S_{12}|\geq 4 \sqrt{2} (b_5+1)e^{1/2}C_2^{1/2} (K_2^2+K_4^2)^{1/2} C_3^{-1}L^{3/2}_n \delta / n )\nonumber\\
 &\leq&8L_n \exp\{-\delta^2/(b_1L^{-1}_nn+b_2\delta)\}+ 6L^2_n\exp \{ -b_6L^{-3}_n n\}.
\label{eq62}
\ed

To bound $S_{13}$, note that
\bg
 |S_{13}|=\ba^T\bU_n({\bU}^{-1}-{\bU}^{-1}_n)\bU\ba
\leq \|{\bU}_n\|^2 \cdot \|{\bU}^{-1}-{\bU}^{-1}_n\|\cdot\|\ba\|_2^2.
\label{eq63}
\ed
Then it follows from Lemmas \ref{lemma6}, Lemma \ref{lemma7}, (\ref{eq58}), (\ref{eq61}), (\ref{eq63}) and the union bound of probability that there exist $b_3$, $b_4$ and $b_6$ such that
 \bg
  &&P(|S_{13}|\geq 4e C_2 (K_2^2+ K_4^2) (b_5+1)^2 C_3^{-2} L^3_n \delta /n)\nonumber\\
 &\leq& 6L^2_n\exp\{-\delta^2/ (b_3L_n^{-1}n + b_4 \delta)\}+ 6L^2_n\exp \{ -b_6L^{-3}_n n\}.
  \label{eq64}
 \ed

Hence, combining (\ref{eq59}), (\ref{eq62}) and (\ref{eq64}), there exist some positive constants $s_1$, $s_2$ and $s_{3}$ such that
\bg
&&P\left( |S_1| \geq s_1 L^2_n \delta^2/n^2+s_2L^{3/2}_n \delta/n+s_3L^{3}_n \delta/n \right)\nonumber\\
&\leq& 16 L_n \exp\{-\delta^2/(b_1 L^{-1}_n n+b_2 \delta)\}+6L^2_n\exp \{ -\delta^2/(b_3 L^{-1}_n  n +b_4 \delta) \} \nonumber\\
&&+18 L^2_n\exp\{-b_6 L^{-3}_n n \}.
 \label{eq65}
\ed
Similarly, we can prove that there exist positive constants $s_4$, $s_5$ and $s_6$ such that
\bg
&&P\left( |S_2| \geq s_4 L^2_n \delta^2/n^2+s_5 L^{3/2}_n \delta/n+s_6 L^{3}_n \delta/n \right)\nonumber\\
&\leq& 8 L_n \exp\{-\delta^2/(b_1 L^{-1}_n n+b_2 \delta)\}+6 L^2_n\exp \{ -\delta^2/(b_3 L^{-1}_n  n +b_4 \delta)\}\nonumber\\
&&+18 L^2_n\exp\{-b_8 L^{-3}_n n \}.
 \label{eq66}
\ed

Let $(s_1+s_4) L^2_n \delta^2/n^2+(s_2+s_5)L^{3/2}_n \delta/n+(s_3+s_6)L^{3}_n \delta/n =c_2L_nn^{-2\kappa}$ for any given $c_2>0$ (e.g., take $\delta= c_2L_n^{-2} n^{1-2\kappa} / (s_3+s_6) $). There exist some positive constants $c_3$ and $c_4$ such that
\bg
  &&P\left(|\hat{u}_{nj} - \tilde{u}_j| \geq c_2L_nn^{-2\kappa}\right)\nonumber\\
  &\leq&(24L_n+12L^2_n)\exp\{-c_3n^{1-4\kappa}L^{-3}_n\}+36L^2_n\exp\{-c_4L^{-3}_n n\}.
\label{eq67}
\ed
Then Theorem \ref{thm1}(i) follows from the union bound of probability.

We now prove part (ii). Note that on the event
$$
\mathcal{A}_n\equiv\left\{\max_{j\in\mathcal{M}_*} |\hat{u}_{nj} - \tilde{u}_j| \leq c_1\xi L_nn^{-2\kappa}/2\right\},
$$
by Proposition \ref{prop1}, we have
\bg
\hat{u}_{nj} \geq c_1\xi L_nn^{-2\kappa}/2,~~~\mbox{for all } j\in\mathcal{M}_*.
\label{eq68}
\ed
Hence, by choosing $\tau_n =c_1\xi L_nn^{-2\kappa}/2$, we have $\mathcal{M}_*\subset\hat{\mathcal{M}}_{\tau_n}$.
On the other hand, by the union bound of probability, there exist positive constants $c_6$ and $c_7$, such that
\bgn
 P(\mathcal{A}_n^c) & \leq& s_n\left\{(24L_n+12L^2_n)\exp(-c_6n^{1-4\kappa}L^{-3}_n)+36L^2_n\exp(-c_7L^{-3}_n n)\right\},
\edn
and Theorem \ref{thm1}(2) follows. $\Box$

 {\bf Proof of Theorem \ref{thm2}}. Let
$$\tilde{\balpha}=\arg\min\limits_{\balpha} \E[(Y-\bQ\balpha)^2],$$
where $\bQ=(\bQ_1,\cdots,\bQ_{p})$ is a  $2pL_n$-dimensional vector of functions.
Then we have $$\E[ \bQ^T(Y- \bQ \tilde{\balpha})]={\bf0}_{2pL_n},$$
where ${\bf 0}_{2pL_n}$ is a $2pL_n$-dimension vector with all entries 0.
This implies
\bgn
 \| \E[\bQ^T Y] \|_2^2 &=& \tilde{\balpha}^T \bSigma^2 \tilde{\balpha}
 \leq \lambda_{\max}( \bSigma) \tilde{\balpha}^T \bSigma \tilde{\balpha},
\edn
recalling $\bSigma = \E[\bQ^T\bQ] $.
It follows from orthogonal decomposition that $\Var(\bQ \tilde{\balpha}) \leq \Var(Y)$ and $\E[\bQ \tilde{\balpha}] = \E[Y]$ (recall the inclusion of the intercept term). Therefore,
\bgn
 \tilde{\balpha}^T \bSigma \tilde{\balpha} \leq \E[Y^2]= O(1),
\edn
and
\bg
 \| \E[\bQ^T Y]\|_2^2=O(\lambda_{\max}(\bSigma)).
 \label{eq69}
\ed
Note that by the definition of $\tilde{u}_j $,
\bgn
 \sum_{j=1}^{p} \tilde{u}_j &=& \sum_{j=1}^p \E[Y\bQ_j] \left((\E[\bQ^T_j \bQ_j])^{-1}  - \left[\begin{array}{cc} (\E[\bB^T\bB])^{-1} & {\bf 0} \\ {\bf 0} & {\bf 0} \end{array} \right] \right) \E[\bQ_j^TY]\\
&\leq&\max_{1\leq j\leq p}\lambda_{\max}\{(\E[\bQ^T_j\bQ_j])^{-1}\}\sum_{j=1}^{p}\|\E[\bQ^T_j Y]\|_2^2\\
 &=&\max_{1\leq j\leq p}\lambda_{\max}\{(\E[\bQ^T_j\bQ_j])^{-1}\}\| \E[\bQ^T Y]\|_2^2.
\edn
By Lemma~\ref{lemma6} and (\ref{eq69}), the last term is of order $O(L_n\lambda_{\max}(\bSigma))$. This implies that the number of $\{j: \tilde{u}_j >\delta L_n n^{-2\kappa}\}$ cannot exceed $O(n^{2\kappa}\lambda_{\max}(\bSigma))$ for any $\delta>0$.

On the set
\bgn
 \mathcal{B}_n=\left\{\max_{1\leq j\leq p}|\hat{u}_{nj} - \tilde{u}_j|\leq \delta L_n n^{-2\kappa}\right\},
\edn
the number of $\{j:\hat{u}_{nj}>2\delta L_n n^{-2\kappa}\}$ cannot exceed the number of $\{j:\tilde{u}_j>\delta L_n n^{-2\kappa}\}$, which is bounded by $O(n^{2\kappa}\lambda_{\max}(\bSigma))$. By taking $\delta=c_5/2$, we have
\bgn
 P\left\{|\hat{\mathcal{M}}_{\tau_n}|\leq O(n^{2\kappa}\lambda_{\max}(\bSigma))\right\}\geq P(\mathcal{B}_n).
\edn
Then the desired result follows from Theorem \ref{thm1}(i). $\Box$


\begin{thebibliography}{}
\bibitem [Antoniadis and Fan(2001)]{Antoniadis2001} Antoniadis, A. and Fan, J. (2001). Regularized wavelet
approximations (with discussion). {\em Jour. Ameri. Statist.
Assoc.}, {\bf 96}, 939-967.

\bibitem[Candes and Tao(2007)]{Candes2007} Candes, E. and Tao, T. (2007). The Dantzig selector:
statistical estimation when $p$ is much larger than $n$ (with
discussion), \ANNALS, {\bf 35}, 2313-2404.
%
\bibitem[Fan, Feng and Song(2011)]{Fan2011} Fan, J., Feng, Y. and Song, R. (2011). Nonparametric
independence screening in sparse ultra-high-dimensional additive
models, \JASA, {\bf 106}, 544-557.
%
%
\bibitem[Fan and Li(2001)]{Fan2001} Fan, J. and Li, R. (2001). Variable selection via
nonconcave penalized likelihood and its oracle properties, \JASA,
{\bf 96}, 1348-1360.
%
\bibitem[Fan and Lv(2008)]{Fan2008} Fan, J. and Lv, J. (2008). Sure independence screening for
ultrahigh dimensional feature space (with discussion), \JRSSB, {\bf
70}, 849-911.
%
%
\bibitem[Fan and Song(2010)]{Fan2010} Fan, J. and Song, R. (2010). Sure independence screening
in generalized linear models with NP-dimensionality, \ANNALS, {\bf
38}, 3567-3604.

\bibitem[Fan, Zhang and Zhang(2001)]{Fan22001} Fan, J., Zhang, C. and Zhang, J. (2001). Generalized
likelihood ratio statistics and wilks phenomenon, \ANNALS, {\bf 29},
153-193.
%
\bibitem[Fan and Zhang(2008)]{Fan22008} Fan, J. and Zhang, W. (2008).   Statistical methods with
varying
        coefficient models. {\em Stat. Interface.}, {\bf 1},
        179-195. 
%
\bibitem[Frank and Friedman(1993)]{Frank1993} Frank, I.E. and Friedman, J.H. (1993). A statistical view
of some chemometrics regression tools (with discussion), {\em
Technometrics}, {\bf 35}, 109-148.
%
\bibitem[Hall and Miller(2009)]{Hall2009} Hall, P. and Miller, H. (2009). Using Generalised Correlation
to Effect Variable Selection in Very High Dimensional Problems,
{\JCGS}, {\bf 18}, 533-550.
%
\bibitem[Hall, Titterington and Xue(2009)]{Hall22009}Hall, P., Titterington, D. M. and Xue, J. H. (2009).
Tilting methods for assessing the influence of components in a
classifier, {\JRSSB}, {\bf 71}, 783-803.
%
\bibitem[Harrison and Rubinfeld(1978)]{Harrison1978} Harrison, D. and Rubinfeld, D. (1978). Hedonic housing
prices and the demand for clean air, {\em J. Environ. Econ. Manag.}
{\bf 5}, 81-102.

\bibitem[Hastie and Tibshirani(1990)]{Hastie} {Hastie, T.} and {Tibshirani, R.} (1990). {\em Generalized
Additive Models}, London: Chapman \& Hall.
%
\bibitem[Hastie and Tibshirani(1993)]{Hastie1993} {Hastie, T.} and {Tibshirani, R.} (1993).
Varying-coefficient models. \JRSSB, {\bf 55}, 757-796.

%
\bibitem [Li, \etal(2012)]{Li2012} Li, G, Peng, H., Zhang, J. and Zhu, L.  (2012).  Robust
rank correlation based screening.  {\em Ann. Statist.}, {\bf 40},
1846-1877.

%
\bibitem[Li, Zhong and Zhu(2012)]{Li22012} Li, R., Zhong, W. and Zhu, L. (2012). Feature screening
via distance correlation learning. \JASA, to appear.

%
\bibitem[Lian(2011)]{Lian2011} Lian, H. (2011). Flexible shrinkage estimation in
high-dimensional varying coefficient models, {\em manuscript}.
%
\bibitem[Stone(1982)]{Stone1982} Stone, C.J. (1982). Optimal glaobal rates of convergence for nonparametric regression, \ANNALS, {\bf 10}, 1040-1053.
%
\bibitem[Tibshirani(1996)]{Tibshirani1996} Tibshirani, R. (1996). Regression shrinkage and selection via the lasso, \JRSSB, {\bf 58}, 267-288.
    %

\bibitem[van der Vaart and Wellner(1996)]{van1996} van der Vaart, A.W. and  Wellner, J.A. (1996).  Weak Convergence and
         Empirical Processes.  Springer, New York.

%
\bibitem[Vershynin(2011)]{Vershynin2011} Vershynin R. (2011). Introduction to the non-asymptotic analysis of random matrices, {\em manuscript}.
%

\bibitem[Wang, Li and Huang(2008)]{Wang2008} Wang, L., Li, H., and Huang, J. Z. (2008).
Variable selection in nonparametric varying-coefficient models for
analysis of repeated measurements. {\JASA}, {\bf 103}, 1556-1569.

\bibitem[Yuan and Lin(2006)]{Yuan2006} Yuan, M. and Lin, Y. (2006). Model selection and estimation in regression with grouped variables. {\em Jour. Roy. Statist.
    Soc. B}, {\bf 68}, 49-67.

\bibitem[Zhang(2010)]{Zhang2010} Zhang, C.-H. (2010), Nearly unbiased variable selection under minimax concave penalty, \ANNALS, {\bf 38}, 894-942.
%
\bibitem[Zhao and Li(2010)]{Zhao2010}Zhao, D. S., and Li, Y. (2010). Principled Sure Independence Screening for
Cox Models With Ultra-High-Dimensional Covariates.  {\em
manuscript}, Harvard University.
%
\bibitem[Zou(2006)]{Zou2006} Zou, H. (2006) The adaptive Lasso and its oracle properties, \JASA, {\bf 101}, 1418-1429.
%
\bibitem[Zou and Hastie(2005)]{Zou2005} Zou, H. and Hastie, T. (2005). Addendum: Regularization and variable selection via the Elastic net, \JRSSB, {\bf 67}, 301-320.
%
\bibitem[Zou and Li(2008)]{Zou2008} Zou, H. and Li, R. (2008). One-step sparse estimates in nonconcave penalized likelihood models, \ANNALS, {\bf 36}, 1509-1533.
\end{thebibliography}
\end{document}